\documentclass[12 pt]{amsart}

\usepackage{amscd, amssymb}
\usepackage{hyperref}
\usepackage{graphicx}
\input xy
\xyoption{all}

\newtheorem{tm}{Theorem}
\newtheorem{pr}[tm]{Proposition}
\newtheorem{lm}[tm]{Lemma}

\theoremstyle{definition}
\newtheorem{dfn}[tm]{Definition}

\theoremstyle{remark}

\newcommand{\proj}{\mathbb P}

\newcommand{\rarr}{\rightarrow}
\newcommand{\oh}{{\mathcal{O}}}
\newcommand{\com}{\mathbb{C}}

\newcommand{\Z}{\mathbb{Z}}
\newcommand{\R}{\mathbb{R}}
\newcommand{\C}{\com}

\newcommand{\K}{\mathcal K}
\newcommand{\E}{\mathcal E}
\renewcommand{\L}{\mathcal{L}}
\newcommand{\rs}{H^0}
\newcommand{\rh}{H^1}
\renewcommand{\O}{\mathcal O}
\newcommand{\p}{\mathfrak p}
\newcommand{\s}{\mathfrak s}
\DeclareMathOperator{\Aut}{Aut}
\DeclareMathOperator{\Def}{Def}
\DeclareMathOperator{\coker}{coker}
\DeclareMathOperator{\im}{Im}
\DeclareMathOperator{\id}{Id}
\DeclareMathOperator{\dist}{dist}
\newcommand{\PP}{\mathbb P}
\newcommand{\open}{disk }

\newcommand{\eqq}{\stackrel{\sim}{=}}

\newcommand{\hb}{\hbar}

\begin{document}

\title{Disk enumeration on the quintic 3-fold}
\author{R. Pandharipande, J. Solomon, and J. Walcher}
\date{May 2007 }

\begin{abstract}
Holomorphic disk invariants with boundary in the real Lagrangian
of a quintic 3-fold are calculated by localization and proven
mirror transforms. A careful discussion of the underlying
virtual intersection theory is included. The generating function
for the disk invariants is shown to satisfy an extension of the
Picard-Fuchs differential equations associated to the mirror
quintic. The Ooguri-Vafa multiple cover formula is used to define
virtually enumerative disk invariants.  The results may also be
viewed as providing a virtual enumeration of real rational curves
on the quintic.
\end{abstract}

\maketitle

\pagestyle{plain}
\setcounter{section}{-1}

\tableofcontents

\section{Introduction}
\subsection{Complex curve enumeration}\label{mm}
Let $Q\subset \com\proj^4$ be a nonsingular quintic hypersurface.
The virtual
count $n_{d}$ of rational algebraic
curves of degree $d>0$ on $Q$ admits a computation
via Gromov-Witten theory and mirror symmetry.

Let $N_{d}$ denote the genus 0  Gromov-Witten invariant
of $Q$ in degree $d$.
The count
$n_{d}$ is
defined by the Aspinwall-Morrison formula \cite{am},
\begin{equation*}
\sum_{d>0}  N_d e^{dT} = \sum_{d>0} \sum_{k>0}
 {n_d}{k^{-3}} e^{kdT}.
\end{equation*}
The connection between
 $n_d$ and actual curve counting
on $Q$ is discussed in \cite{pan}.

The mirror symmetry
prediction
of Candelas, de la Ossa, Green, and
Parkes \cite{cogp} relates the genus 0 potential
$$\mathcal{F}(T)= \frac{5}{6}T^3 + \sum_{d>0} N_d e^{dT}$$
to hypergeometric series.
Let $I_i(t)$ be defined by
\begin{equation}\label{ggt}
\sum_{i=0}^3 I_i H^i = \sum_{d=0}^{\infty} e^{(H+d)t}
\frac{\Pi_{r=1}^{5d} (5H+r)}{\Pi_{r=1}^d(H+r)^5} \ \ \text{mod} \
H^4.
\end{equation}
The functions $I_i(t)$ are a basis of solutions of the
Picard-Fuchs differential equation
\begin{equation} \label{pfq}
\Big( \frac{d}{dt}\Big)^4 I - 5 e^t \Big(5 \frac{d}{dt}+1\Big)
\Big(5 \frac{d}{dt}+2\Big) \Big(5 \frac{d}{dt}+3\Big) \Big(5
\frac{d}{dt}+4\Big) I = 0.
\end{equation}
Let the variables $T$ and $t$ be related by $T(t)= I_1/I_0 \ (t)$.
The prediction,
\begin{equation*}
\mathcal{F}(T(t))= \frac{5}{2} \Big( \frac{I_1}{I_0}(t)
\frac{I_2}{I_0}(t) - \frac{I_3}{I_0}(t) \Big),\end{equation*}
was later proven via localization on the space of
genus 0 stable maps to $\com\proj^4$ \cite{g1, ko,lly}.

\subsection{Disk enumeration}
Let $Q\subset \com\proj^4$ be a nonsingular quintic hypersurface
{defined
over} $\mathbb{R}$.
Let $\omega$ be the symplectic form on $Q$ obtained
from the Fubini-Study metric.
Complex conjugation determines an anti-holomorphic
involution on $Q$ with
fixed locus equal to the set of real points $Q_{\mathbb{R}}$.
The inclusion
$$Q_{\mathbb{R}} \subset Q$$
is  Lagrangian with respect to $\omega$.

The enumeration of disks  with boundaries in Lagrangian
submanifolds plays a basic role in open string theory and
has been studied mathematically in several contexts.
The subject is not a direct extension of the theory
of stable
maps. New issues such as orientation play a crucial role.
We provide
here a complete calculation
of the disk invariants of $Q$ with
boundary in the real Lagrangian $Q_{\mathbb{R}}$.

An early treatment of disk enumeration occurs in the
 construction of the Fukaya category
\cite{FOOO}. Disk enumeration is required to define the
differentials of the Floer complex. However, a symplectic
invariant via disk enumeration is not defined in \cite{FOOO}.
Only
the cohomology of the Floer complex is invariant.

Symplectic disk invariants have been defined with respect to
the real Lagrangian associated to
an anti-holomorphic involution in \cite{jake1}.
A previous definition in the presence of a torus action preserving
the real Lagrangian (not directly applicable to $Q$)
can be found in \cite{katzliu,melissa}.
We will follow here the definitions of \cite{jake1}.

Let $N_d^{disk}$ for $d$ odd denote the degree $d$ disk invariant of $Q$
with boundary in $Q_{\mathbb{R}}$. For a discussion of
even degree, see Sections \ref{ssec:ms} and \ref{ssec:even}. In fact,
$N_d^{disk}$
depends on a choice of a $Spin$ structure on
$Q_\R.$ However, changing the $Spin$ structure only effects the sign of
$N_d^{disk}$ uniformly for all $d$. 
Our conventions are fixed by choosing $N_1^{disk}$ to be
positive.  Let $\mathcal{F}^{disk}$ denote the disk potential,
$$\mathcal{F}^{disk}(T)= \sum_{d \ odd}   N^{disk}_d e^{dT/2}.$$
Our main result is a calculation of $\mathcal{F}^{disk}$. Define
\[
J(t) = 2 \sum_{d \ odd}\ e^{dt/2}\ \frac{(5d)!!}{(d!!)^5}.
\]
The function $J(t)$ is a solution of the Picard-Fuchs equation \eqref{pfq} with
an added inhomogenous term,
\[
\Big( \frac{d}{dt}\Big)^4 J - 5 e^t
\Big(5 \frac{d}{dt}+1\Big)\!\Big(5 \frac{d}{dt}+2\Big)\!\Big(5
\frac{d}{dt}+3\Big)\!\Big(5
\frac{d}{dt}+4\Big) J = \frac{15}{8} e^{t/2}.
\]
Alternatively, $J(t)$ may be obtained by evaluation at $H=1/2$ 
of the (non-truncated) hypergeometric solution \eqref{ggt}
of the homogeneous Picard-Fuchs equation,
$$J(t)=  30 \Bigg[ \sum_{d=0}^{\infty} e^{(H+d)t}
\frac{\Pi_{r=1}^{5d} (5H+r)}{\Pi_{r=1}^d(H+r)^5} \Bigg]_{H=\frac{1}{2}}\ .
$$

\begin{tm} \label{tm1} Via the mirror map $T(t)= I_1/I_0(t)$,
\begin{equation*}
\mathcal{F}^{disk}(T(t))= \frac{J(t)}{I_0(t)}.
\end{equation*}
\end{tm}

The disk invariants $N_d^{disk}$ are typically fractional.
Following the strategy of curve enumeration, 
virtual disk counts $n_d^{disk}$ are defined by the
Ooguri-Vafa formula.
\begin{dfn} \label{tm2} We define the counts $n_d^{disk}$ by
\begin{equation*}\sum_{d\ odd}  N^{disk}_d e^{dT/2} =
\sum_{d\  odd}\ \sum_{k\ odd}
 {n^{disk}_{d/k}}{k^{-2}} e^{kdT/2}.
\end{equation*}
\end{dfn}
\noindent 
Definition \ref{tm2} is justified by the multiple cover calculation of
Proposition \ref{der} in Section \ref{mcf}. 
The contribution of
$k$-fold covers of a disk in the appropriate
local Calabi-Yau geometry is 
$k^{-2}$.
We conjecture the invariants $n_d^{disk}$ to
be integers.

\subsection{Real curve enumeration}
Each holomorphic disk with boundary in $Q_{\mathbb{R}}$
can be reflected by the Schwartz principle
to yield a real rational curve in $Q$. Conversely,
real curves mapping to $Q$
 of odd degree may be halved to yield two disks \cite{jake1}.
The virtual number of real rational curves of odd degree in $Q$ may be
defined by
$$n_d^{real} = \frac{1}{2} n_d^{disk}.$$
Again, $n_d^{real}$ vanishes for $d$ even.
A table of values can be found at the end of the paper.

\subsection{Mirror Symmetry}\label{ssec:ms}
Let $\L$ be a $U(1)$ bundle with flat connection $A_0$ over the
Lagrangian submanifold $Q_\R \subset Q.$ The triple 
$$\O = (Q_\R,\L,A_0)$$
determines an object of the Fukaya category of $Q.$ Homological mirror symmetry
\cite{hams} predicts 
the existence of a corresponding object $\O^\vee$
of the derived category of coherent sheaves on the mirror quintic $Q^\vee.$  
The holomorphic Chern-Simons functional of
$\O^\vee$ is predicted to be mirror to the standard Chern-Simons functional of
$\O$ with corrections from disk instantons \cite{W}. In the following, 
we briefly explain
how the mirror correspondence between functionals leads to an enumerative
correspondence.

Assume for simplicity that $\O^\vee$ is a holomorphic vector bundle. Denote
the underlying complex vector bundle by $V$, 
and let $A_{\O^\vee}$ be the connection on
$V$ defining the holomorphic structure of $\O^\vee.$ Let $t$ denote the complex
moduli parameter of $Q^\vee,$ and let $\Omega_t$ denote the holomorphic 
3-form determined by $t.$ The holomorphic Chern-Simons functional 
of $V$ depends
on a second complex connection on $V,$ which we denote $A.$ We view $A$ as
a connection 1-form relative to $A_{\O^\vee}.$ 
Define the holomorphic Chern-Simons
functional by
\[
L^\vee(A,t)  = \int_{Q^\vee} \text{Tr}\left(A \wedge 
\bar\partial_{A_{\O^\vee}} A 
+ \frac{2}{3} A \wedge A \wedge A \right) \wedge \Omega_t.
\]
Critical points of $L^\vee$ are holomorphic connections ---
complex connections with vanishing $(0,2)$ component of their curvature.
$L^\vee$ is constant on connected components of orbits of the complex gauge
group of $V.$

For the Chern-Simons functional with instanton corrections, we will
require the following terminology for holonomy.
If $$E\rightarrow B$$
 is a bundle
with connection $\theta,$ and $P$ is a parametrized path in $B$, we denote the
holonomy of $\theta$ around $P$ by $\text{Hol}(P,\theta).$ If $B \subset Q$ and 
\[
f:(D^2,\partial D^2) \rightarrow (Q,B),
\]
we write 
\[
P_f = \{ f|_{\partial D^2} : \partial D^2 \rightarrow B \}.
\]
Moreover, if $f$ is holomorphic, we define $\eta_f$ to be the sign of $f$
coming from the determinant line of the Cauchy-Riemann operator. 

In defining the Chern-Simons functional with instanton corrections, we will
not try to be entirely precise, but rather give an intuitive picture. 
Let $M
\subset Q$ be a totally real submanifold homologous to $Q_\R$,
 and let $C^4$ be a
chain satisfying 
$\partial C^4 = M - Q_\R$.
Let 
$$\L_M \rightarrow M$$ be a $U(1)$
bundle, and let $A_M$ be a connection on $\L_M.$ 
Let $T$ denote the complexified Kahler moduli parameter of $Q$ and let
$\omega_T$ denote the associated complexified Kahler form. Define the
Chern-Simons function with instanton corrections by
\begin{align}
L(M,A_M,T) & = \int_{M} A_M \wedge (d A_M - \omega) - \int_{Q_\R} A_0 \wedge d
A_0 + \int_{C^4} \omega_T \wedge \omega_T \notag\\
& \quad + \sum_{\substack{f: (D^2,\partial
D^2) \rightarrow (Q,M) \\ \bar \partial f = 0}} \eta_f \text{Hol}(P_f, A_M)
\exp\left(\int_{D^2} f^* \omega_T \right) \notag\\
& \quad - \sum_{\substack{f: (D^2,\partial
D^2) \rightarrow (Q,Q_\R) \\ \bar \partial f = 0}} \eta_f \text{Hol}(P_f, A_0)
\exp\left(\int_{D^2} f^* \omega_T\right) \notag\\
\label{eq:cs} & \quad + \sum_{\substack{f: S^2 \rightarrow Q \\ \bar\partial f =
0,\, z \in f^{-1}(C^4)}} \eta_f \exp\left(\int_{S^2} f^* \omega_T\right) 
\end{align}
A Lagrangian submanifold $M$ with vanishing obstruction chains in the sense
of \cite{FOOO} and a flat connection $A_M$ together constitute a critical point
of $L.$ $L$ is constant on orbits of the Hamiltonian symplectomorphism
group. The corrections from closed instantons intersecting $C^4$ are necessary
to compensate for codimension one bubbling where all the energy of a disk
instanton is transferred to a sphere bubble \cite{S}.

$L$ and $L^\vee$ are Lagrangians defining a pair of dual quantum field
theories \cite{W}. Their critical values are physically significant and should
be topological invariants.
Via the mirror transformation, which expresses $T$ as a function of $t,$ 
the value of $L$ at a critical point should be calculable
from the value of $L^\vee$ at a mirror critical point. However, we must find
non-trivial choices of critical points $A$ and $(M,\L_M,A_M)$ 
corresponding under mirror
symmetry for all values of $T = T(t)$.

In our case,
such $(M,\L_M,A_M)$ can be found using the geometry of
the anti-holomorphic involution.
We choose 
$$M = Q_\R.$$
Since $H_1(Q_\R) = \Z/2,$, we can choose
$(\L_M,A_M)$ to be the flat $U(1)$ bundle with monodromy opposite to
$(\L,A_0).$ Since $Q_\R$ is the fixed point set of an anti-holomorphic
involution, $Q_\R$ is a critical point of $L$ with either flat
bundle. The mirror to such a choice of $(\L_M,A_M)$ should be a unique up to
gauge transformation holomorphic connection $A$ on $V$ not gauge equivalent to
$A_{\O^\vee}.$ In \cite{Walcher}, a heuristic argument is given to show
$L^\vee(A,t)$ is in fact given by $J(t)$. Even after specifying which
holomorphic structure $A$ induces, $L^\vee(A,t)$ is only defined up to a period
of $\Omega_t.$ Indeed, changing the choice of $A$ by a complex gauge
transformation not isotopic to the identity changes $L^\vee(A,t)$ by a period of
$\Omega_t.$ Therefore, $L^\vee(A,t)$ should satisfy an extension of the
Picard-Fuchs equation for $Q^\vee.$ Similarly, changing the choice of $C^4$ in
the definition of $L$ changes $L$ by a multiple of the first derivative of
$\mathcal F(T).$

From the preceding discussion, we see both $L$ and $L^\vee$ are
essentially relative functionals, depending either on a pair of connections
$A,\,A_{\O^\vee},$ or a pair of totally real submanifolds with $U(1)$ bundle,
$(M,\L_M)$ and $(Q_\R,\L).$ Therefore only disks of odd degree, for which
the difference in monodromy of $A_M$ and $A_0$ cancels the negative sign in 
definition \eqref{eq:cs}, contribute to the physically significant critical
value of $L^\vee.$ The contributions of even degree disks cancel due to this
sign. A priori, some other physical value may depend on the
even degree disks. However, on mathematical grounds, 
even degree disks appear not to lead to interesting invariants, see Section
\ref{ssec:even}.

\subsection{Past and future work}
The first number $n_1^{disk}$ was calculated in \cite{jake1}.
Theorem 1 was predicted in \cite{Walcher} via low degree
graphs sums and string heuristics.
Our technique of proof uses the fully equivariant mirror
correspondence of Givental \cite{g1}.
A previous
application can be found in \cite{grzas} where
disk enumeration for (noncompact) local geometries
was considered.
The Ooguri-Vafa \cite{oovafa} multiple cover formula of Definition \ref{tm2} 
is by now established
in many settings, see \cite{katzliu,melissa}.

We have chosen the quintic 3-fold
as a first case of study, but the methods of the
paper are applicable much more generally. It will be 
interesting to see which aspects of the solution persist.

\subsection{Acknowledgements}
We thank P. Biran, T. Graber, D. Kazhdan, N. Nekrasov, L. Polterovich, P.
Seidel, E. Shustin, G. Tian, R. Vakil, E. Zaslow, and A. Zinger for related
conversations.

J.~S.~ would like to thank E. Farjoun and the Hebrew University of
Jerusalem for their warm hospitality during the preparation of
the paper.
R.~P.~ was partially supported by a Packard foundation fellowship
and NSF grant DMS-0500187.
J.~S.~ was partially supported by NSF grant DMS-0111298.
J.~W. was partially 
supported by NSF grant PHY-0503584.

\section{Disk invariants}\label{di}
\subsection{Overview}\label{dj}
We recall the definition of the disk invariant $N_d^{disk}$
of the quintic from \cite{jake1}.  
Our conventions for conjugation, real structures, and
stable disk maps are discussed in Sections \ref{conn}.
The Euler class approach to $N_d^{disk}$ is presented in
Section 
\ref{ssec:q3d}.

\subsection{Definitions}
Fix a symplectic manifold $(X,\omega)$ of real
dimension less than or equal to $6$ with an anti-symplectic
involution $\phi,$
\[
\phi^*\omega = - \omega.
\]
The fixed points $L = Fix(\phi)$ define a Lagrangian submanifold
of $X.$ A $Pin$ structure and, if $L$ is orientable, an
orientation on $L$ induce a natural relative orientation on the
moduli space $\overline{M}_{D}(X/L,\beta)$ of stable \open 
maps to
$(X,L)$ of degree $\beta$. 
Since $\overline{M}_{D}(X/L,\beta)$ is an orbifold with
corners, the definition of cohomology
classes on the moduli space yielding an analog of Gromov-Witten
theory is not immediately clear. 
However, using $\phi,$ 
certain corners of $\overline{M}_{D}(X/L,\beta)$ may
be eliminated.

More precisely, the boundary  of 
$\overline{M}_{D}(X/L,\beta)$ consists generically
of stable \open  maps with two \open
components. Replacing {\em one} of the two components 
by the image under $\phi$ yields
another two component map. We define an equivalence relation $\sim$ on
the boundary of $\overline{M}_{D}(X/L,\beta)$ based on 
this correspondence. For
certain components of the boundary, the relation $\sim$ preserves
orientation. After quotienting by $\sim$ on these components, 
we obtain a new moduli
space $\widetilde{M}_{D}(X/L,\beta)$ with fewer corners which is still
relatively orientable. On $\widetilde{M}_{D}(X/L,\beta)$, 
many interesting cohomology classes can be defined. 
Consequently, a set
of invariants are obtained
of the triple $(X,\omega,\phi)$  reminiscent of
standard Gromov-Witten invariants in many respects.

In good
situations, the invariants obtained
from $\widetilde{M}_{D}(X/L,d)$
 are actually enumerative. For example, 
Welschinger's signed counts of real curves \cite{We1,We2}
arise as specializations of the theory \cite{jake1}.

\subsection{Conventions}\label{conn}
\subsubsection{Coordinates}
Let $z_0, \ldots, z_4$ be homogeneous coordinates on $\com\proj^4$.
The standard complex conjugation $c_{\id}$ on $\com\proj^4$ is
$$[z_0, z_1,z_2,z_3,z_4] \stackrel{c_{\id}}{\mapsto} [\overline{z}_0,
\overline{z}_1,\overline{z}_2,\overline{z}_3,\overline{z}_4].$$
Each $g\in \mathbb{PGL}_5$ yields an anti-holomorphic
involution
$$c_g=g^{-1} \circ c_{\id} \circ g : \com\proj^4 \rarr \com\proj^4$$
equivalent to $c_{\id}$.
In particular, the anti-holomorphic involution $c$,
\begin{equation}\label{ccc}
[z_0,z_1,z_2,z_3,z_4] \stackrel{c}{\mapsto} [\overline{z}_0,
\overline{z}_2,\overline{z}_1,\overline{z}_4,\overline{z}_3],
\end{equation}
is associated to the matrix
$$\left( \begin{array}{ccccc}
1 & 0 & 0 & 0 & 0 \\
0 & 1 & 1 & 0 & 0 \\
0 & i & -i & 0 & 0 \\
0 & 0 & 0 & 1& 1 \\
0 & 0 & 0 & i & -i \\ \end{array} \right)\ .$$ Let
$\com\proj^4_{\mathbb R} \subset \com\proj^4$ denote the fixed
points of $c$. The involution $c$ will be most convenient for our
calculation of disk invariants.

\subsubsection{Real geometry}
A homogeneous polynomial $F(z_0,z_1,z_2,z_3,z_4)$ on $\com\proj^4$
is defined over $\mathbb{R}$ if
$$\overline{F(z)} = F(c(z)).$$
For example,
\begin{equation}\label{frz}
z_1+z_2
\ \ \text{and}\ \
iz_1-iz_2
\end{equation}
are both real linear polynomials.

A subvariety of $V\subset \com\proj^4$ is defined over $\mathbb{R}$
if the ideal $I(V)$ is generated by real homogeneous functions.
The lines
$$L=\{\ [0,z_1,z_2, 0 ,0]\ |\ z_1,z_2\in \com\ \},$$
$$
L'=\{\ [0,0,0,z_3,z_4]\ |\ z_3,z_4\in \com \ \}$$
are both defined over $\mathbb{R}$.

The involution $c$ lifts canonically to the line bundles
\[
\mathcal{O}_{\com\proj^4}(k) \rarr \com\proj^4.
\]
The linear polynomials \eqref{frz} are elements of
\[
H^0(\com\proj^4,\mathcal{O}_{\com\proj^4}(1))_{\mathbb{R}} \subset
H^0(\com\proj^4,\mathcal{O}_{\com\proj^4}(1))_{\mathbb{C}},
\]
the space of real sections.

\subsubsection{Maps}
\label{ccc2}
Let $u,v$ be homogeneous coordinates on $\com\proj^1$. Let
$$c: \com\proj^1 \rarr \com\proj^1$$
be the anti-holomorphic involution defined by
\begin{equation*}
[u, v] \stackrel{c}{\mapsto} [\overline{v},\overline{u}].
\end{equation*}
The $c$-fixed points, $\com\proj^1_{\mathbb R}\subset \com\proj^1$ form
 a circle.

A holomorphic disk map
$$f:(D, \partial D) \rarr (\com\proj^4, \com\proj^4_{\mathbb R})$$
can be reflected by the Schwartz principle to yield an algebraic
map
$$\tilde{f}: \com\proj^1 \rarr \com\proj^4.$$
By definition, the degree $d$ of the disk map equals the
degree of $\tilde{f}$.

The map $\tilde{f}$ satisfies the following real condition
\begin{equation}\label{ned}
\tilde{f}\circ c = c \circ \tilde{f}.
\end{equation}
Conversely, every algebraic map
$$\tilde{f}: \com\proj^1 \rarr \com\proj^4$$
satisfying \eqref{ned} yields two disk maps with
boundary $\partial D$ equal to $\com\proj^1_{\mathbb{R}}$.
The image of $\tilde{f}$ is a real subcurve of $\com\proj^4$.

Similarly, a stable holomorphic disk map $f$ reflects to
a stable genus 0 map $\tilde{f}$ satisfying the real
condition \eqref{ned} with respect to the natural
extension of $c$ to degenerations of $\com\proj^1$.
In fact, stability for $f$ can be defined by stability
for $\tilde{f}$.
We will use the notation
$$f: (D, \partial D) \rarr (\com\proj^4, \com\proj^4_{\mathbb R})$$
also for the stable case where $D$ and $\partial D$ are possibly
reducible. However, $\partial D$ is always connected.

In the odd degree case, {\em every} stable genus 0
map to $\com\proj^4$ defined over $\mathbb{R}$
is obtained by reflection.

\subsubsection{Moduli}
Let $M_D(\com\proj^4/\com\proj^4_{\mathbb R},d)$ denote the moduli
space of unpointed disk maps of {odd} degree $d$. Reflection
yields an \'etale double cover of smooth orbifolds
$$\epsilon: M_D(\com\proj^4/ \com\proj^4_{\mathbb R},d)
\rarr M_{\mathbb R}(\com\proj^4,d)
$$
where $M_{\mathbb R}(\com\proj^4,d)$ denotes the moduli space of
unpointed genus $0$ algebraic maps defined over $\mathbb{R}$. The
real dimension of $M_{\mathbb{R}}(\com\proj^4,d)$ is $5d+1$. In
fact, $\epsilon$ is an orientation double cover \cite{jake1}.

Let $\overline{M}_D(\com\proj^4/\com\proj^4_{\mathbb R},d)$ denote
the compactification of the moduli space $M_D(\com\proj^4/\com\proj^4_{\mathbb
R},d)$ by stable disk maps,
 and let
$\overline{M}_{\mathbb{R}}(\com\proj^4,d)$ denote the space of
unpointed genus $0$ algebraic stable maps defined over
$\mathbb{R}.$ The moduli space  
$\overline{M}_D(\com\proj^4/\com\proj^4_{\mathbb R},d)$ is a
smooth orbifold with corners. In fact,
$\epsilon$ extends to finite smooth map
\[
\bar\epsilon : \overline{M}_D(\com\proj^4/\com\proj^4_{\mathbb
R},d) \rarr \overline{M}_\R(\com\proj^4,d),
\]
mapping the corners of
$\overline{M}_D(\com\proj^4/\com\proj^4_{\mathbb R},d)$ to the
boundary divisor of
$\overline{M}_\R(\com\proj^4,d).$ The cardinality of the fiber
over a real stable map with $n_o$ components fixed by $c$ is
$2^{n_o}.$

In Section \ref{dj}, the construction of the closed orbifold
\[
\widetilde M_D(\com\proj^4/\com\proj^4_\R,d) = 
\overline{M}_D(\com\proj^4/\com\proj^4_{\mathbb R},d) / \sim
\]
was outlined.
A detailed argument is given in Section \ref{sec:ecf} in the proof
of Proposition \ref{pr:co}. The equivalence relation $\sim$ identifies
the corners of $\overline{M}_D(\com\proj^4/\com\proj^4_{\mathbb R},d)$ in such
a way that 
the map $\bar \epsilon$ descends to an \'etale double cover
\[
\tilde\epsilon : \widetilde M_D(\com\proj^4/\com\proj^4_\R,d) \rarr
\overline{M}_\R(\com\proj^4,d).
\]
There is a natural inclusion
$$\overline{M}_{\mathbb R}(\com\proj^4,d)\subset
\overline{M}_{\com}(\com\proj^4,d)$$ in the space of unpointed
stable genus 0 algebraic maps defined over $\com$. The real
dimension of $\overline{M}_{\mathbb{R}}(\com\proj^4,d)$ is $5d+1$.

\subsection{Euler class formula} \label{ssec:q3d}
Let $Q \subset \C \proj^4$ be a nonsingular quintic 
hypersurface defined over $\R$
with symplectic form obtained from the Fubini-Study metric.
 An anti-symplectic involution 
$$\phi: Q \rightarrow Q$$ is defined by complex
conjugation. The Lagrangian $Fix(\phi)$ is the real locus
$Q_\R.$ 

We consider maps from the holomorphic disk
 $D$ to $Q$ of odd degree with boundary lying in $Q_\R$,
$$\overline{M}_D(Q/Q_\R,d) \subset 
\overline{M}_D(\com\proj^4/\com\proj^4_{\mathbb
R},d).$$
Since the expected dimension of the
moduli space of maps to $Q/Q_R$ is 0, the relevant Gromov-Witten
invariant $N_d^{disk}$ is simply the virtual 
cardinality.


In order 
to calculate $N_d^{disk},$ following \cite{jake1}, we first reformulate
$N_d^{disk}$ as the integral of an Euler class of an
obstruction bundle over the moduli space
$\widetilde{M}_D(\com\proj^4/\com\proj^4_{\mathbb R},d).$ Such
integrals may be studied via fixed point localization. A similar
approach was used by Kontsevich \cite{ko} in the closed case.

Let $\hat F_d$ be the real vector bundle over
$\overline{M}_D(\com\proj^4/\com\proj^4_{\mathbb R},d)$ with fiber
$$\hat F_d|_{[f:(D,\partial D)\rarr(\com\proj^4,\com\proj^4_{\mathbb{R}})]}
= H^0(C,\tilde{f}^*\mathcal{O}_{\com\proj^4}(5))_{\mathbb{R}}$$
where
$$[\tilde{f}:C\rarr \com\proj^4] \in \overline{M}_{\mathbb R}(\com\proj^4,d)$$
is the stable rational map obtained from the stable disk map via
reflection, and $H^0(C,\tilde{f}^*\mathcal{O}_{\com\proj^4}(5))_{\mathbb{R}}$
denotes real sections.

The vector bundle $\hat F_d$ is of real rank $5d+1$ and is
oriented on $\overline{M}_D(\com\proj^4/\com\proj^4_{\mathbb
R},d)$ by Lemma 8.7 of \cite{jake1}. The integral of the
Euler class $e(\hat F_d)$ over
$\overline{M}_D(\com\proj^4/\com\proj^4_{\mathbb R},d)$ is not
well defined because the space
$\overline{M}_D(\com\proj^4/\com\proj^4_{\mathbb R},d)$ has
non-empty boundary, and $\hat F_d$ is not trivial near the
boundary. However,
$\hat F_d$ naturally descends to a
vector bundle $$F_d \rightarrow
\widetilde{M}_D(\com\proj^4/\com\proj^4_{\mathbb R},d).$$
 Neither
$F_d$ nor $\widetilde{M}_D(\com\proj^4/\com\proj^4_{\mathbb R},d)$
are orientable. Let $\L$ denote the local system defined
by the determinant of the tangent bundle of the moduli space
$\widetilde{M}_D(\com\proj^4/\com\proj^4_{\mathbb R},d).$ In Lemma \ref{lm:ro}
of
Section \ref{sec:ecf}, we prove 
\[
\det{F_d} \simeq \L
\]
as topological bundles.
A $Spin$ structure on $Q_\R$ determines the
choice of the isomorphism uniquely up to scaling by a positive
constant. Hence, the Euler class
\[
e(F_d) \in
H^{5d+1}(\widetilde{M}_D(\com\proj^4/\com\proj^4_{\mathbb R},d),
\L)
\]
is well-defined. Since
$\widetilde{M}_D(\com\proj^4/\com\proj^4_{\mathbb R},d)$ is a
closed orbifold, the integral
$$\int_{\widetilde{M}_D(\com\proj^4/\com\proj^4_{\mathbb R},d)}
e(F_d)$$ is well-defined. In Section \ref{sec:ecf}, we obtain the
following result.
\begin{tm}\label{tm3}
For $d$ odd,
$$N_{d}^{disk} =
\int_{\widetilde{M}_D(\com\proj^4/\com\proj^4_{\mathbb R},d)}
e(F_d).$$
\end{tm}
We prove Theorem \ref{tm3} using the symplectic
virtual moduli cycle. The same technique can be used to prove the analogous
well-known result for the closed invariants.

\subsection{\texorpdfstring{$N_d^{disk}$}{Invariant} in even degree}
\label{ssec:even}
A stable disk map of even degree may still be reflected to obtain
an even degree real genus 0 stable map. However, not all stable
genus 0 maps of even degree defined over $\mathbb{R}$ are so
obtained. Stable maps defined over $\mathbb{R}$ with domains
having no real points cannot be halved.

The even disk invariant $N_d^{disk}$ is not well-defined
without the addition of the contributions of real curves
without real points. If such contributions were incorporated,
$N_d^{disk}$ would be expressible as the Euler class of an odd dimensional real
bundle and hence would presumably vanish. Hence, the definition $N_d^{disk} = 0$
for $d$ even.

\subsection{Dependence on \texorpdfstring{$Q$}{Q}}
Our formula for $N_d^{disk}$ is {\em independent} of the
quintic $Q\subset \com\proj^4$ defined over $\mathbb{R}$.
Since the calculation is done on $\com\proj^4$, some
information is possibly lost.
More precisely, let 
$$\epsilon:H_1(Q_{\mathbb{R}},\mathbb{Z}/2\mathbb{Z}) \rightarrow 
H_1(\com\proj^4_{\mathbb{R}}, \mathbb{Z}/2\mathbb{Z}) \eqq \mathbb{Z}/2\mathbb{Z}.$$
The invariant $N_d^{disk}$ is an integral over all stable
disk maps 
$$f:(D,\partial D) \rarr (Q, Q_{\mathbb{R}})$$
of degree $d$ with boundary $\partial D$ determining a
class in 
$$\epsilon^{-1}(\overline{1})\subset H_1(Q_{\mathbb{R}},
\mathbb{Z}/2\mathbb{Z}).$$
If $\epsilon$ is an isomorphism, as is the case, for example,
for the Fermat quintic
$$Q=(z_0^5+z_1^5+z_2^5+z_3^5+z_4^5),$$
then there is no loss of information.
If, however, $\epsilon$ has a kernel, more refined disk invariants
of $(Q,Q_{\mathbb{R}})$ may sometimes be  defined for
$$\gamma\in \epsilon^{-1}(\overline{1})\subset 
H_1(Q_{\mathbb{R}}, \mathbb{Z}/2\mathbb{Z}).$$ In the latter
case, 
$$N_d^{disk} = \sum_{\gamma \in \epsilon^{-1}(\overline{1})} 
N_{d,\gamma}^{disk}.$$

\section{Torus actions}

\subsection{Tori}
Let $\mathbf{T}$ denote the complex numbers of unit modulus,
$$\mathbf{T}= \{ \xi\in \com  \ | \ |\xi|=1\ \}.$$
The torus $\mathbf{T}^5$ acts
diagonally $\com^5$. A
$\mathbf{T}^5$-action on $\com\proj^4$ is obtained by projectivization,
and canonical lifts
to the line bundles
$$\mathcal{O}_{\com\proj^4}(k) \rarr \com\proj^4$$
are obtained.
There is a canonically induced
translation action of $\mathbf{T}^5$ on
$\overline{M}_{\com}(\com\proj^4,d)$.

Let $\zeta_i \in \com\proj^4$ denote the $\mathbf{T}^5$-fixed
points,
$$\zeta_0=[1,0,0,0,0],\ \zeta_1=[0,1,0,0,0],\ \ldots,\ \zeta_5=[0,0,0,0,1].$$
The involution $c$ fixes $\zeta_0$ and permutes the
others,
$$\zeta_1 \stackrel{c}{\leftrightarrow}\zeta_2,\ \
\zeta_3 \stackrel{c}{\leftrightarrow}\zeta_4.$$
Hence, $\zeta_0$ is the unique real $\mathbf{T}^5$-fixed point.

Consider the rank
2 subtorus
$\mathbf{T}^2 \subset \mathbf{T}$
acting by
$$(\xi_1,\xi_2) \cdot [z_0,z_1,z_2,z_3,z_4] =
[z_0, \xi_1 z_1, \overline{\xi}_1 z_2 , \xi_2 z_3,
\overline{\xi}_2 z_4].$$
Since $\mathbf{T}^2$ preserves $\com\proj^4_{\mathbb R}$,
translation defines a $\mathbf{T}^2$-action on the
moduli spaces
$\overline{M}_D(\com\proj^4/ \com\proj^4_{\mathbb R},d)$ and
$\widetilde{M}_D(\com\proj^4/ \com\proj^4_{\mathbb R},d)$.

The algebraic torus $(\com^*)^2$ acts on $\com\proj^4$
by complexifying the action of $\mathbf{T}^2$,
$$(\xi_1,\xi_2) \cdot [z_0,z_1,z_2,z_3,z_4] =
[z_0, \xi_1 z_1, \xi^{-1}_1 z_2 , \xi_2 z_3,
\xi^{-1}_2 z_4].$$
Of course, $(\com^*)^2$ acts on $\overline{M}_{\com}(\com\proj^4,d)$
by translation.

\subsection{Equivariant weights}
We follow the equivariant weight conventions of \cite{g1,pan} for the
torus
$\mathbf{T}^5$.

Let $\lambda_i$ be the $\mathbf{T}^5$-equivariant cohomology class
determined by Chern class of the restriction of
$\mathcal{O}_{\com\proj^4}(1)$ to $\zeta_i$,
$$\lambda_i=c_1(\mathcal{O}_{\com\proj^4}(1)_{\zeta_i}) \in
H^*_{\mathbf{T}^5}(\text{pt}).$$
The classes $\lambda_i$ generate,
$$H^*_{\mathbf{T}^5}(\text{pt})= \mathbb{Q}[\lambda_0,\ldots,\lambda_4].$$
The tangent weights of $\com \proj^4$
at the point $\zeta_i$ are
$\{ \lambda_i-\lambda_j  \}_{j\neq i}$.

Let $\lambda, \lambda'$ be the
generators of $H^*_{\mathbf{T}^2}(\text{pt})$
defined by the pull-back
$$\rho^*: H^*_{\mathbf{T}^5}(\text{pt})
\rarr H^*_{\mathbf{T}^2}(\text{pt})$$
and the equations
$$
\rho^*(\lambda_1)=-\rho^*(\lambda_2) = \lambda,\
\rho^*(\lambda_3)=-\rho^*(\lambda_4)=\lambda'.$$
The pull-back $\rho^*(\lambda_0)$ vanishes.
For notational convenience, we will often omit the
pull-back $\rho^*$ and write
\begin{equation}\label{laz}
\lambda_0=0, \
\lambda_1=-\lambda_2= \lambda, \ \lambda_3=-\lambda_4=\lambda'.
\end{equation}

\subsection{Localization}
The genus 0 Gromov-Witten invariants $N_d$ have been
calculated in \cite{g1,ko,lly}  via
localization on $\overline{M}_{\com}(\com\proj^4,d)$
with respect to the $\mathbf{T}^5$-action. We
will calculate  $N_d^{disk}$  via localization
on $\widetilde{M}_D(\com\proj^4/ \com\proj^4_{\mathbb R},d)$
with respect to the
$\mathbf{T}^2$-action.

\section{Localization calculation of \texorpdfstring{$\mathcal{F}^{disk}$}{disk potential}}

\subsection{Overview}
Let $d$ be odd. The $\mathbf{T}^2$-action on the moduli space
${\widetilde{M}_D(\com\proj^4/\com\proj^4_{\mathbb R},d)}$ lifts
canonically to the vector bundle $F_d$. We calculate the integral
$$N_{d}^{disk} =
\int_{\widetilde{M}_D(\com\proj^4/\com\proj^4_{\mathbb R},d)}
e(F_d)$$ by localization with respect to the
$\mathbf{T}^2$-action.

The localization calculation is similar in flavor to
the genus 0 Gromov-Witten calculation of $Q$ in \cite{ko}.
However, two new issues arise:
\begin{enumerate}
\item[(i)] The $\mathbf{T}^2$-action has fixed loci in
${\widetilde{M}_D(\com\proj^4/\com\proj^4_{\mathbb R},d)}$ with
moving images in $\com\proj^4$. \item[(ii)] The equivariant
restriction of $e(F_d)$ to the $\mathbf{T}^2$-fixed locus depends
upon the orientation of $F_d$.
\end{enumerate}
Issue (i) is handled by identifying the non-rigid
contributions with the equivariant correlators
$S_Q$
studied
by
Givental \cite{g1,pan}.
Issue (ii) requires an explicit evaluation of the
signs occurring in the orientation. The derivation
is presented in Section \ref{signs}.

The sum over $\mathbf{T}^2$-fixed point loci
required for the localization formula is executed
in two steps. Subsums with fixed intersection type
with $\com\proj^4_{\mathbb{R}}$ are evaluated
using Givental's equivariant mirror transformation
for $S_Q$. Finally, the sum over intersection types
is evaluated explicitly after appropriate equivariant
specialization. The interaction of the orientation signs
with the localization sum is an interesting
aspect of the calculation.
The outcome is a proof of Theorem 1.

\subsection{\texorpdfstring{$\mathbf{T}^2$}{Torus}-fixed disk maps}
We first study the $\mathbf{T}^2$-fixed locus of the moduli space
of stable disk maps.
Let
$$[f: (D,\partial D)
\rarr
(\com\proj^4,\com\proj^4_{\mathbb R})]
\in{\overline{M}_D(\com\proj^4/\com\proj^4_{\mathbb R},
d)^{\mathbf{T}^2}}$$
be a $\mathbf{T}^2$-fixed map.

The boundary
$\partial D$ distinguishes a minimal, $c$-invariant, central
curve $P\subset C$ of the domain of the reflected map
$$[\tilde{f}:
C \rarr \com\proj^4] \in \overline{M}_{\mathbb{R}}(\com\proj^4,d)$$
satisfying
$\partial D = P_{\mathbb{R}}$.
The {\em central degree}  of $f$ is the
degree of the restriction
$$\tilde{f}_{P}: P \rarr \com\proj^4.$$
The {central degree} $p$ is positive, odd, and
bounded by $d$.
The moduli point
$$[\tilde{f}_P]\in \overline{M}_\com(\com\proj^4,p)$$
is fixed for the full complexified
action of
$(\com^*)^2$ on $\com\proj^4$.

\begin{lm} The two lines $L,L'\subset \com\proj^4$
 are the only $(\com^*)^2$-invariant
curves of odd degree defined over $\mathbb R$
in $\com\proj^4$.
\end{lm}

\begin{proof} A real $(\com^*)^2$-invariant subcurve must
lie in one of the two planes
$$\{\ [z_0,z_1,z_2,0,0] \ | \ z_0,z_1,z_2 \in \com \ \},$$
$$\{\ [z_0,0,0,z_3,z_4] \ | \ z_0,z_3,z_4 \in \com \ \}.$$
In the first case, $L$ is the only real $(\com^*)^2$-invariant
line. Moreover, all non-linear $(\com^*)^2$-orbits are
of degree 2.
The argument in the second case is identical.
\end{proof}

A node of $P_{\mathbb R}$ must map via $\tilde{f}$ to the unique
real fixed point $\zeta_0\in \com\proj^4$. Since $\tilde{f}(P)$
must equal either $L$ or $L'$, $\tilde{f}(P)$ can not contain
$\zeta_0$. Hence, $P_{\mathbb R}$ cannot contain a node. We obtain
the following result.

\begin{lm}
The central curve $P$ is $\com\proj^1$ and
$$\tilde{f}_P:P \rarr L \ \text{or}\ L'$$
is a Galois cover of odd degree $p$.
\end{lm}

The original disk map $f$ is obtained from one half of
$\tilde{f}$. Hence one half of $P$ is selected by $D$.
A half of $P$ determines a pair $(\zeta,p)$
where
$$\zeta\in \{\zeta_1, \zeta_2,\zeta_3,\zeta_4 \}$$
is a non-real fixed point and
$p$ is the central degree.

The data $(\zeta,p)$ is the termed the {\em intersection type}
of $f$ with the real Lagrangian $\com\proj^4_{\mathbb{R}} \subset
\com\proj^4$. The half of $P$ selected by $D$ is the
{\em intersection disk}.

While we have analyzed 
$\overline{M}_D(\com\proj^4/\com\proj^4_{\mathbb R},d)^{\mathbf{T}^2}$, 
we
are actually interested in 
$\widetilde{M}_D(\com\proj^4/\com\proj^4_{\mathbb R},d)^{\mathbf{T}^2}$.
If fact, we have proven
$$\overline{M}_D(\com\proj^4/\com\proj^4_{\mathbb R},d)^{\mathbf{T}^2}=
\widetilde{M}_D(\com\proj^4/\com\proj^4_{\mathbb R},d)^{\mathbf{T}^2}$$
since the $\mathbf{T}^2$-fixed maps are {\em not} corner points
of $\overline{M}_D(\com\proj^4/\com\proj^4_{\mathbb R},d)$.

\subsection{Intersection disk term}
The localization calculation of
$$N_{d}^{disk} =
\int_{\widetilde{M}_D(\com\proj^4/\com\proj^4_{\mathbb R},d)}
e(F_d)$$ is sum over the contributions of the $\mathbf{T}^2$-fixed
loci. We may separate the contributions by intersection type,
$$N_d^{disk} = \sum_{i=1}^4 \sum_{p \ odd } \text{Cont}_{(\zeta_i,p)}
(N_d^{disk}).$$
The {\em intersection disk term} $I(\zeta_i,p)$ of
$\text{Cont}_{(\zeta_i,p)}
(N_d^{disk})$
is the contribution of the unique
$\mathbf{T}^2$-fixed map
$$f: (D,\partial D) \rarr (\com\proj^4,\com\proj^4_{\mathbb{R}})$$
incident to $\zeta_i$ with central degree $p$ and domain
consisting only of the intersection disk.

Define the rational function $C_p(\lambda,\lambda')$ of
degree $0$ by the following
formula,
\begin{eqnarray*}
C_p(\lambda,\lambda')
& = & \frac{(-1)^\frac{p-1}{2}
}{p}
\frac{2\lambda}{p}
\frac{\frac{(5p)!!}{p! p!!}
(\frac{\lambda}{2p})^{p}}
{ \prod_{i=0}^{(p-1)/2} ( (1-\frac{2i}{p})\lambda-\lambda')
( (1-\frac{2i}{p})\lambda+\lambda')
 }.
\end{eqnarray*}

\begin{lm}\label{ft1}
For an appropriate choice of $Spin$ structure on $Q_\R$, we have
\[
I(\zeta_1,p)=I(\zeta_2,p)=C_p(\lambda,\lambda'), \qquad
I(\zeta_3,p)=I(\zeta_4,p)=C_p(\lambda',\lambda).
\]
Changing the $Spin$ structure changes the formulas by $-1$ for all
$p.$
\end{lm}

The proof of Lemma \ref{ft1} is given in Section \ref{signs}. The
most interesting aspect is the calculation of the prefactor
$(-1)^{\frac{p-1}{2}}$ obtained from the orientations of the
moduli space $\overline{M}_D(\com\proj^4/
\com\proj^4_{\mathbb{R}},p)$ and the bundle $F_p$.

\subsection{Givental's correlator \texorpdfstring{$S_Q$}{}}
Let $\overline{M}_{0,2}(\com\proj^4,r)$ be the moduli space
of 2-pointed stable complex algebraic maps to
$\com\proj^4$ of genus 0 and degree $r$. Let
$$e_i: \overline{M}_{0,2}(\com\proj^4,r) \rarr \com\proj^4$$
be the evaluation at the $i^{th}$ marking, and
let $\psi_i$ denote the Chern class of the $i^{th}$
cotangent line.
Let
$$E_r \rarr  \overline{M}_{0,2}(\com\proj^4,r)$$
be the complex vector bundle with
fiber
$$E_r|_{[f:C\rarr \com\proj^4 ]}
= H^0(C,\tilde{f}^*\mathcal{O}_{\com\proj^4}(5))_{\com}.$$

Following the notation of Section 2.2 of \cite{pan},
Givental's equivariant correlator $S_Q$ for the
torus $\mathbf{T}^2$
is defined by
\begin{equation}\label{ert}
S_Q(T,\hbar)=
\frac{1}{5H}
\sum_{r\geq 0} e^{(H/\hbar+r) T}   e_{2*}(
\frac {c_{\text{top}}(E_r)}{\hb-\psi_2})\in H^*_{\mathbf{T}^2}(\com\proj^4)
\end{equation}
where $H$ is the hyperplane class,
$$H= c_1(\mathcal{O}_{\com\proj^4}(1)) \in H^*_{\mathbf{T}^2}(\com\proj^4).$$
The sum in \eqref{ert} is over {\em all} non-negative integers $r$.
The unstable degree 0 term is defined by
\begin{equation*}
\frac{1}{5H}
e_{2*}( \frac{c_{\text{top}}(E_0)}{\hb-\psi_2}) =1.
\end{equation*}

Let $[\zeta_i] \in H^*_{\mathbf{T}^2}(\com\proj^4)$
denote the Poincar\'e dual of the class of the fixed point.
For classes $\mu,\nu \in H^*_{\mathbf{T}^5}(\com\proj^4)$,
let
$$\langle \mu, \nu \rangle \in \mathbb{Q}[\lambda, \lambda']$$
denote the equivariant intersection pairing. For example,
$$\langle H, [\zeta_i]\rangle = \lambda_i$$
following convention \eqref{laz}.

The intersection pair of the equivariant correlator will
arise in the localization analysis:
\begin{eqnarray*}
\langle S_Q(T,\hbar), [\zeta_i] \rangle &  = &
\frac{1}{5\lambda_i}
\sum_{r\geq 0} e^{(\lambda_i/\hbar+r) T}
\int_{\overline{M}_{0,2}(\com\proj^4,r)}
\frac {c_{\text{top}}(E_r)}{\hb-\psi_2} e_2^*([\zeta_i])  \\
& = & \frac{\hb^{-1}}{5\lambda_i}
\sum_{r\geq 0} e^{(\lambda_i/\hbar+r) T}
\int_{\overline{M}_{0,1}(\com\proj^4,r)}
\frac {c_{\text{top}}(E_r)}{\hb-\psi_1} e_1^*([\zeta_i])
\end{eqnarray*}
where the string equation is used in the second line.
In degree $0$, the unstable 1-pointed term is defined
by the second equality.

\subsection{Contributions of type \texorpdfstring{$(\zeta_i,p)$}{(zeta i,p)}}
The $\mathbf{T}^2$-fixed loci of the moduli space
$\widetilde
{M}_D(\com\proj^4/\com\proj^4_{\mathbb{R}},d)$
of type $(\zeta_i,p)$ may be quite complicated.
However, every map
$$[f]\in \widetilde{M}_D
(\com\proj^4/\com\proj^4_{\mathbb{R}},d)^{\mathbf{T}^2}$$
of type $(\zeta_i,p)$
ends in the same intersection disk.
By expanding the localization formula,
the intersection disk term $I(\zeta_i,p)$
can be factored out of
$\text{Cont}_{(\zeta_i,p)}(N^{disk}_d)$
by removing the intersection disk from $f$.

What remains after the intersection disk is removed from $f$?
In fact, {\em every} genus 0 stable complex map
$$[f']\in e_1^{-1}(\zeta_i) \subset
\overline{M}_{0,1}(\com\proj^4,r)^{(\com^*)^2}$$
can be found.
The stable disk map $f$ is obtained by attaching
the $(\zeta_i,p)$-intersection disk to $f'$
at the marking.

A direct unraveling of the localization formulas
yields the following fundamental result.
Let
$$\text{Cont}_{(\zeta_i,p)}(\mathcal{F}^{disk}) =
\sum_{d \ odd} e^{dT/2}
\text{Cont}_{(\zeta_i,p)}(N^{disk}_d).$$

\begin{lm} We have
$$\text{\em  Cont}_{(\zeta_i,p)}(\mathcal{F}^{disk}) =
\langle S_Q(T,\frac{2}{p}\lambda_i), [\zeta_i] \rangle
\cdot I(\zeta_i,p)$$
for $1\leq i \leq 4.$
\end{lm}

\begin{proof}
The sum on left side can be indexed more conveniently as
$$\text{ Cont}_{(\zeta_i,p)}(\mathcal{F}^{disk}) =
\sum_{r\geq 0} e^{(\frac{p}{2}+r) T} 
\text{Cont}_{(\zeta_i,p)}(N^{disk}_{{p+2r}}).$$
The right side of the equality may be expanded as
$$
\sum_{r\geq 0} e^{(\frac{p}{2}+r) T}
\int_{\overline{M}_{0,1}(\com\proj^4,r)}
\frac {c_{\text{top}}(E_r)}{\frac{2}{p}\lambda_i -\psi_1} e_1^*([\zeta_i]) 
\cdot \frac{I(\zeta_i,p)}{(5\lambda_i)(\frac{p}{2} \lambda_i)}.
$$
The Lemma is obtained from the equality
\begin{equation}
\label{jlem}
\text{Cont}_{(\zeta_i,p)}(N^{disk}_{{p+2r}}) =
\int_{\overline{M}_{0,1}(\com\proj^4,r)}
\frac {c_{\text{top}}(E_r)}{\frac{2}{p}\lambda_i -\psi_1} e_1^*([\zeta_i]) 
\cdot \frac{I(\zeta_i,p)}{(5\lambda_i)(\frac{p}{2} \lambda_i)}.
\end{equation}

To prove \eqref{jlem}, we apply $\mathbf{T}^2$-localization to
the integral on the right.
We do {\em not} fully expand the $\mathbf{T}^2$-localization
formula. As was discussed previously, understanding the
geometry of the individual $\mathbf{T}^2$-fixed
loci is difficult as there are positive dimensional
families of $\mathbf{T}^2$-fixed curves.
However, \eqref{jlem} has a much simpler proof. Since both
sides are expressed as $\mathbf{T}^2$-residue integrals
by localization, we need only match the geometries.

First, the 
$\mathbf{T}^2$-fixed loci of the two sides of \eqref{jlem}
are in bijective correspondence.
Given a $\mathbf{T}^2$-fixed locus on the right, the addition of the
intersection disk $I(\zeta_i,p)$ at the marking $1$ 
yields a $\mathbf{T}^2$-fixed locus of the right side. The reverse
direction is obtained by stripping the intersection disk.

Second, since the $\mathbf{T}^2$-fixed loci on the left and right
are both nonsingular, the correspondence induces an {\em isomorphism} 
of $\mathbf{T}^2$-fixed loci up to
the automorphism factor of the intersection disk.

Finally, we must match 
the $\mathbf{T}^2$-fixed obstruction theories. Let
$$[f_D:D= C \cup I(\zeta_i,p)\rightarrow \com\proj^4]\in 
\widetilde{M}_D
(\com\proj^4/\com\proj^4_{\mathbb{R}},p+2r)^{\mathbf{T}^2}
$$
be a map with 
$$[f_C:C \rarr \com\proj^4] \in \text{ev}_1^{-1}(\zeta_i)\subset
\overline{M}_{0,1}
(\com\proj^4,r)^{\mathbf{T}^2}.$$
By the normalization sequence,
$$0 \rarr  F_{p+2r}|_{[f_D]} \rarr E_{r}|_{[f_C]} \oplus
F_{p}|_{[I(\zeta_i,p)]} \rarr \oh_{\com\proj^4}(5)|_{\zeta_i} \rarr 0.$$
Hence, the numerator in the residue integral on the left of
\eqref{jlem} is
\begin{equation}\label{jyy}
e(F_{p+2r}) 
= c_{\text {top}}(E_r) \cdot \frac{e(F_{p})}{5\lambda_i}.
\end{equation}
Similarly, the denominator of the residue integral on the left
of \eqref{jlem} is 
\begin{equation}\label{tyy}
\frac{1}{e(\text{Nor}_{[f_D]})} = \frac{1}{c_{\text{top}}
(\text{Nor}_{[f_C]})}
\frac{c_{\text{top}}(\text{Tan}_{\zeta_i})}
{(\frac{2}{p} \lambda_i -\psi_1) (\frac{2}{p}\lambda_i)}
 \frac{1}
{e(\text{Nor}_{[I(\zeta_i,p)]})}.
\end{equation}
The middle terms are obtained from tangent bundle, node smoothing,
and automorphism factors. Putting \eqref{jyy} and \eqref{tyy}
together, we obtain the exact matching needed for \eqref{jlem}.

The
factorization of
\eqref{jlem} properly reflects the orientation on the
moduli space $\widetilde{M}_D(\com\proj^4/\com\proj^4_{\mathbb{R}},d)$. 
The orientation factorization
is easily obtained from \cite{WW}. 
\end{proof}

Such arguments form the
geometric basis of \cite{g1,pan}.
Though the $\mathbf{T}^2$-action on
$\widetilde{M}_D(\com\proj^4/\com\proj^4_{\mathbb{R}},d)$
has fixed loci corresponding to moving maps, the issue is
completely avoided by the Lemma.

\subsection{Mirror transforms}
We review the mirror transforms relating $S_Q(T,\hb)$
to
$$S^*_Q(t,\hb)= \frac{1}{5H}
\sum_{r\geq 0}
 e^{(\frac{H}{\hb}+r)t}
\frac{
\Pi_{s=0}^{5r}( 5H+s \hb)}{
\Pi_{j=0}^{4}
\Pi_{s=1}^r(H-\lambda_j+s\hb)}$$
following Section 4.4 of \cite{pan}.

The mirror map $T(t)= I_1/I_0(t)$ discussed in Section \ref{mm}
can be written explicitly.
Let
$$F(q)= \sum_{r=0}^\infty q^r \frac{(5r)!}{(r!)^{5}},
\ \
G_l(q)= \sum_{r=1}^{\infty} q^r \frac{(5r)!}{(r!)^{5}}
\Big( \sum_{s=1}^{lr} \frac{1}{s}\Big).$$
Then
\begin{equation*}
T= t+ \frac{5(G_{5}(e^t)-G_1(e^t))}{F(e^t)}
\end{equation*}
is the mirror map.
Exponentiating yields
\begin{equation*}
\exp(T)= \exp(t)\cdot\text{exp}
\left(\frac{5(G_5(e^t)-G_1(e^t))}{F(e^t)}\right).
\end{equation*}

The equivariant mirror transformation for the torus
$\mathbf{T}^2$ is
$$S_Q(T(t),\hb) = \frac{1}{F(e^t)} S_Q^*(t,\hb).$$
Transforms for the equivariant pairings are a direct
consequence,
$$\langle S_Q(T(t),\frac{2}{p} \lambda_i), [\zeta_i]\rangle
= \frac{1}{F(e^t)} \langle S_Q^*(t,\frac{2}{p}\lambda_i),
[\zeta_i]\rangle.$$

\subsection{Theorem 1}
We now complete the calculation of $\mathcal{F}^{disk}$,
\begin{eqnarray*}
\mathcal{F}^{disk}  & = & \sum_{i=1}^4 \sum_{p\ odd} \text{Cont}_{(\zeta_i,p)}
(\mathcal{F}^{disk}) \\
&= & \sum_{i=1}^4 \sum_{p\ odd}
\langle S_Q(T,\frac{2}{p}\lambda_i), [\zeta_i] \rangle
\cdot I(\zeta_i,p) \\
&= & \sum_{i=1}^4 \sum_{p\ odd}\frac{1}{F(e^{t})}
\langle S^*_Q(t,\frac{2}{p}\lambda_i), [\zeta_i] \rangle
\cdot I(\zeta_i,p)\, .
\end{eqnarray*}
The $i=1,2$ summands of the last line together yield
\begin{equation*}
\frac{4}{F(e^t)}\sum_{r\geq 0} \sum_{p \ odd}
e^{(\frac{p}{2}+r)t}
\frac{(-1)^{\frac{p-1}{2}}}{r!(p+r)!p}\frac{
2^{-(p+2r)}  \frac{(5p+10r)!!}{ (p+2r)!!}}
{\prod_{1\leq i\ odd \leq p+2r} (i-px) (i+px)}
\end{equation*}
written in terms of the
homogeneous variable
$$x= \frac{\lambda'}{\lambda}.$$
Similarly, the $i=3,4$ summands together yield
$$\frac{4}{F(e^t)}\sum_{r\geq 0} \sum_{p \ odd}
e^{(\frac{p}{2}+r)t}
\frac{(-1)^{\frac{p-1}{2}}}{r!(p+r)!p}\frac{
2^{-(p+2r)}  \frac{(5p+10r)!!}{ (p+2r)!!}}
{\prod_{1\leq i\ odd \leq p+2r} (i-px^{-1}) (i+px^{-1})}.$$

The final step is to observe the localization calculation of
$\mathcal{F}^{disk}$ is a {\em weight independent} global
integral. Hence, we may evaluate the summation in the
$x\rarr 0$ limit. Only the $i=1,2$ terms
survive
the limit.
$$
\mathcal{F}^{disk} =
\frac{4}{F(e^t)}\sum_{r\geq 0} \sum_{p \ odd}
e^{(\frac{p}{2}+r)t}
\frac{(-1)^{\frac{p-1}{2}}}{r!(p+r)!p}\frac{
2^{-(p+2r)}  \frac{(5p+10r)!!}{ (p+2r)!!}}
{(p+2r)!! (p+2r)!!}
$$
The identity for odd $d$,
\begin{equation}\label{cccz}
\sum_{1\leq p \ odd \leq d}  \frac{(-1)^{\frac{p-1}{2}}}
{\left(\frac{d-p}{2}\right)!(\frac{d-p}{2}+p)!p} =
\frac{2^{d-1}}{(d!!)^2},
\end{equation}
restated in Lemma \ref{idd} below
concludes the proof of Theorem 1,
\begin{equation*}
\mathcal{F}^{disk}(T)= \frac{2}{F(e^t)}\sum_{d \ odd}
e^{{dt}/{2}}
  \frac{(5d)!!}{(d!!)^5}.
\end{equation*}

After regrouping the factors and reindexing the sum,
the identity \eqref{cccz} is equivalent to the following
result.

\begin{lm}\label{idd} For $d$ odd,
$$\sum_{k=0}^d \binom{d}{k} (-1)^k \frac{d}{d-2k} =
(-1)^{\frac{d-1}{2}} \frac{2^{2d-1}}{ \binom{d-1}{\frac{d-1}{2}}}.$$
\end{lm}

\begin{proof} 
An elementary derivation is left to the reader. A geometric
proof is obtained from the multiple cover calculations in 
Section \ref{mcf}.
\end{proof}

\section{Intersection disk terms}\label{signs}
\subsection{Overview}
    We now derive the signs needed in the localization calculation and prove
Lemma \ref{ft1}. The orientations of the moduli spaces and vector
bundles used to calculate $N_d^{disk}$ arise from the natural
orientation of the determinant of the Cauchy-Riemann $Pin$ boundary
value problem developed in \cite{jake1}. Briefly, a Cauchy-Riemann
boundary value problem consists of topological complex vector
bundle $E$ over a Riemann surface with boundary $\Sigma,$ a
totally real sub-bundle $F$ over the boundary $\partial\Sigma,$
and a generalized Cauchy-Riemann operator $d''$ on $E.$ Given a
$Pin$ structure on $F,$ and a choice of orientation on $F$ if 
 orientable, one may define a canonical orientation of the
determinant
\begin{equation*}
\det(d'') := \Lambda^{\mbox{\tiny max}}(\ker d'')\otimes
\Lambda^{\mbox{\tiny max}}(\coker d'')^*.
\end{equation*}
Reversing the $Pin$ structure on $F$ reverses the
canonical orientation \cite[Lemma 2.10]{jake1}.

For the calculations of weights below, we will only be concerned
with the situation where $E \rightarrow D$ is the restriction of
an algebraic vector bundle $\tilde E \rightarrow \C \proj^1$
defined over $\R$. We take $F = \tilde E_\R,$ and
we take $d''$ to be the restriction to $D$ of the $\bar\partial$
operator defined by the holomorphic structure on $\tilde E.$ The
identifications
\[
\ker(d'') = H^0(\C \PP^1,\tilde E)_\R, \qquad \coker(d'') = H^1(\C
\PP^1,\tilde E)_\R
\]
are easily obtained.
So, an orientation of $\det(d'')$ gives an orientation of the
virtual vector space
\[
H^0(\C \PP^1,\tilde E)_\R - H^1(\C \PP^1,\tilde E)_\R.
\]
Note, however, the orientation depends on the choice of 
$$D
\subset \C \PP^1\setminus \C \PP^1_\R.$$
 So, we cannot entirely forget
the origins of our orientation in a boundary value problem. For
convenience, we introduce the notation
\[
\rs(D,\tilde E) := \ker(d''), \qquad \rh(D,\tilde E) :=
\coker(d'').
\]
Section \ref{sec:ws} interprets the symplectic geometric
definition of the orientation algebraically. Section
\ref{sec:wcalc} calculates the localization contributions from a
fixed point of the torus action using the combinatorics of exact
sequences and the formula for the tensor product of real
representations of $S^1.$ The algebraically inclined reader may
safely skip all of Section \ref{sec:ws} besides the statement of
Lemma \ref{lem:aor}.

\subsection{Weights of sections of a line bundle}\label{sec:ws}
Let $\mathbf{T}^n$ denote the real $n$-dimensional torus and let
$\mathfrak{t}^n$ denote its Lie algebra. A weight is a homomorphism of real
vector spaces from $\mathfrak{t}^n$ to $\C.$ Let $V$ be a two-dimensional real
irreducible oriented representation of $\mathbf{T}^n$ and let $\rho$ be the
associated homomorphism
\[
\rho : \mathbf{T}^n \rightarrow \mathbf{GL}(V).
\]  
Let $h : V \rightarrow \C$ be an orientation preserving real linear
homomorphism such that associated homorphism 
\[
\tilde h : \mathbf{GL}(V) \rightarrow \mathbf{GL}(\C)
\]
satisfies
\[
\im(\tilde h\circ\rho) \subset \text{Aut}_\C(\C) \simeq \C^{\times}.
\]
Note that $h$ is defined by these conditions up to homothety and hence $\tilde
h$ is unique. Differentiating $\tilde h\circ \rho$ we obtain the weight of $V.$

We denote by $V_\lambda$ the 2-dimensional real 
oriented representation of $\mathbf{T}^n$ of weight 
$\lambda,$ where $\lambda$ may be fractional. If $\lambda=0$ we denote
by $V_\lambda$ the trivial two dimensional representation. Suppose $V$ is
an oriented real representation of $\mathbf{T}^n$ isomorphic to $V_\lambda.$
A priori, $V$ consists of two data: a $\mathbf{T}^n$ action on the vector space
$V$ and an orientation of the vector space $V.$
However, except in the case $\lambda = 0,$ the $\mathbf{T}^n$ action on $V$
and the knowledge that $V \simeq V_\lambda$ determine the orientation of
$V$. Indeed, if $\lambda \neq 0,$ there is a unique up to
homothety $\mathbf{T}^n$ equivariant isomorphism 
\[
i: V \stackrel{\sim}{\rightarrow} V_\lambda,
\]
The orientation on $V$ must agree with
the orientation induced by $i.$

Let $m$ be a positive odd integer. From \cite{jake1}, in order to define the
canonical orientation of $\det(d'')$ mentioned above, a $Pin$ structure
$\p_{-1}$ on $\O(-1)_\R \rightarrow \partial D$ must be chosen. We fix such a
$Pin$ structure. Set $W_\lambda = V_\lambda \otimes \C.$ Choosing a connected
component of the complement of
\[
\com \proj_\R^1 \subset \C \proj^1 = \C \proj(W_\lambda).
\]
is equivalent to choosing an orientation on $V_\lambda.$ The
action of $\mathbf{T}^n$ on $W_\lambda$ naturally induces an action on
$\rs(D,\O(m))$ when we think of $D$ as the disk that induces the orientation of
$V_\lambda.$ Indeed, as not oriented vector spaces, clearly
\begin{equation}\label{ccff}
\rs(D,\O(m)) \simeq Sym^m(V_\lambda^*).
\end{equation}
We only consider the $m$ odd case.

The main goal of this section is to prove the following Lemma
which examines when isomorphism \eqref{ccff}
preserves
orientation. Equip the vector space
$\rs(D,\O(m))$ with the $\mathbf{T}^n$-action induced
from $V_{-\lambda}.$ Let $\p$ be a $Pin$ structure on $\O(m)_\R
\rightarrow \partial D.$ Since $m$ is odd, $\O(m)_\R
\simeq \O(-1)_\R$ as real topological vector bundles.
\begin{lm}\label{lem:aor}
Assume that $\p$ agrees with $\p_{-1}.$  With respect to the
canonical orientation induced by $\p,$
\begin{equation}\label{eq:ws}
\rs(D,\O(m)) \simeq \bigoplus_{i = 0}^{(m-1)/2} V_{(2i +
1)\lambda}.
\end{equation}
\end{lm}
\begin{proof}
First, the canonical orientation induced by $\p$ can be
expressed as the complex orientation induced by an explicit
complex structure on $\rs(D,\O(m)).$ Indeed, by gluing sections,
we have an exact sequence
\begin{align}\label{eq:glue}
&0 \rightarrow \rs(D,\O(m)) \overset{g^{-1}}{\rightarrow} \notag\\
&\qquad \qquad \rs(D,\O(-1)) \oplus H^0(\C \proj^1,\O((m+1)/2))
\overset{h}{\rightarrow} \C \rightarrow 0.
\end{align}
By definition, after equipping the latter two terms of the
sequence with the complex orientation, the sequence induces the
desired orientation on the first term. Note that $g^{-1},$ the
inverse of the gluing map, is not canonical, but the set of all
choices is connected. So, the induced orientation is well-defined.
To calculate the orientation, we may fix a particular
choice of the gluing map and calculate the induced complex
structure on $\rs(D,\O(m)).$

We digress for a moment to explain the relationship between
different possible descriptions of $\rs(D,\O(m)).$ Let $w_0,w_1$
be standard linear coordinates on $W_\lambda$ such that $w_0,w_1$
are real precisely on the real locus of $W_\lambda.$ Then,
$\rs(D,\O(m))$ may be identified with the vector space of
homogeneous polynomials
\[
p(w_0,w_1) = \sum_{i = 0}^{m} a_i w_0^{d-j} w_1^j, \qquad a_i \in
\R.
\]
Now, we choose new coordinates $z_0,z_1,$ such that
\begin{gather*}
z_0 = w_0 - i w_1, \qquad z_1 = w_0 + i w_1, \\
w_0 = \frac{z_0 + z_1}{2}, \qquad w_1 = i\frac{z_0 - z_1}{2}.
\end{gather*}
Then, since $p(z_0,z_1) = \overline{p(\bar z_1,\bar z_0)},$ we
have
\[
p(z_0,z_1) = \sum_{j = 0}^{m} b_i z_0^{d-j}z_1^j, \qquad b_j =
\bar b_{d-j}.
\]
On the other hand, by definition, $\rs(D,\O(m))$ is the space of
the solutions $\xi$ of the Cauchy Riemann equations on unit disk
\[
D = \{|z| \leq 1 \} \subset \C
\]
satisfying the totally real boundary conditions
\begin{equation*}
\xi(z) \in \R z^{m/2}, \qquad |z| = 1.
\end{equation*}
In the future, we refer to these boundary conditions as $L(m).$
Expanding in power series about $z = 0,$ it is not hard to see
that the boundary conditions imply
\begin{equation}\label{eq:xi}
\xi(z) = \sum_{j = 0}^m b_i z^i, \qquad b_j = \bar b_{d-j}.
\end{equation}
So, solutions $\xi$ arise naturally from homogeneous
polynomials $p(z_0,z_1)$ by trivializing $\O(m)$ over $\{z_0 \neq
0\}$ by the section $z_0^m$ and making the identification $z =
z_1/z_0.$

For the following argument, we view
sections $\xi \in \rs(D,\O(m))$ as in \eqref{eq:xi}. Referring to
exact sequence \eqref{eq:glue}, we construct a gluing map
\[
g : \ker(h) \rightarrow \rs(D,\O(m))
\]
as follows. Suppose we identify the point at which we glue
sections with $\infty \in \C \proj^1.$ Then, by taking the
standard trivialization of the sheaf $\O((m+1)/2)$ over the standard
coordinate chart on $\C \proj^1$ centered at $0,$ we may identify
$\ker(h)$ with the set of polynomials $q(z)$ of degree less than
or equal to $(m-1)/2.$ Let $$\beta : D \rightarrow \R$$ be a cutoff
function depending only on $|z|$ such that $\beta(0) = 1$ and
$\beta(z) = 0$ for $|z| > 1/2.$ Given a polynomial $q,$ we define
a pre-gluing
\[
\tilde q(z) = \beta(z) q(z).
\]
Let
\[
P : C^\infty((D,\partial D),(\C,L(m))) \rightarrow \rs(D,\O(m))
\]
denote the $L^2$ projection. We define $g$ by
\[
g(q) = P \tilde q.
\]
Due to our choice of gluing map, it is easy to calculate that
\[
g(b z^j) = \bar b z^{m-j} + b z^j.
\]
Let $\xi$ be as in \eqref{eq:xi} and let $I$ denote the complex
structure $g$ induces on $\rs(D,\O(m)).$ It follows that
\[
(I\xi)(z) = i\sum_{j = 0}^{(m-1)/2} b_j z^j - b_{m-j}z^{m-j}.
\]

Finally, we compute the weights of $\rs(D,\O(m)).$ Since, by
definition of the action on $\rs(D,\O(m)),$ we have identified $D$
with the hemisphere of $\C \proj^1$ inducing the intrinsic
orientation of $V_\lambda,$ we must have that $z$ transforms by $z
\mapsto e^{-2\lambda} z.$

Since $z = \frac{z_1}{z_0},$ it follows that in coordinates
$z_0,z_1,$ the action of $\mathbf{T}^n$ takes the form
\[
z_0 \mapsto e^{i\lambda}z_0, \qquad z_1 \mapsto e^{-i\lambda} z_1.
\]
In other words, $\mathbf{T}^n$ acts by $e^{I(m-2j)\lambda}$ on the section
\[
z_0^j z_1^{m-j} + z_0^{m-j} z_1^j \in \rs(D,\O(m)).
\]
The decomposition $\eqref{eq:ws}$ follows.
\end{proof}

\subsection{Localization contribution of an isolated fixed point}
\label{sec:wcalc}

Let $X$ and $Y$ be two oriented real vector spaces, and let $e_i
\in X$ and $f_j \in Y$ be oriented bases. There are different ways
the tensor product $X \otimes Y$ can inherent an orientation from
$X$ and $Y.$ We  use the right-to-left
lexicographical ordering convention. We take
\[
e_1\otimes f_1,\ e_2\otimes f_1, \ldots, e_1 \otimes f_2,\ e_2
\otimes f_2, \ldots
\]
as an oriented basis of $X\otimes Y.$

If $X$ is even dimensional then the orientation of $X
\otimes Y$ is independent of the orientation of $Y,$ and vice
versa.

As before, let $V_\lambda$ be the
2-dimensional real representation of $\mathbf{T}^n$ with weight
$\lambda.$ The following result is a straightforward linear algebra
calculation.
\begin{lm}
With respect to the right-to-left lexicographical orientation of
the tensor product,
\begin{equation}\label{eq:t}
V_\alpha \otimes V_\beta = V_{\alpha + \beta} \oplus
V_{\alpha-\beta}.
\end{equation}
\end{lm}

\subsection {Proof of Lemma \ref{ft1}}
 Let 
$$[f : (D,\partial D) \rightarrow
(\C \proj^4, \com \proj_\R^4)]\in
\overline{M}_D(\com\proj^4/\com\proj^4_{\mathbb R},d)$$
 denote a $\mathbf{T}^2$-fixed disk of type
$(\zeta_1,p).$ After permuting indices, the proof given below
applies to the other possible intersection types as well.

First, we calculate the equivariant Euler class of  the tangent space to
$\overline{M}_D(\com\proj^4/\com\proj^4_{\mathbb R},d)$ at $[f]$
which we denote $N_f.$ We use the deformation exact sequence
\begin{equation}\label{5t3}
0 \rightarrow \Aut(D) \rightarrow \Def(f) \rightarrow \Def(D,f)
\rightarrow 0.
\end{equation}
Here, $\Def(f)$ denotes the space of first-order deformations of
the map $f,$ and $\Def(f,D)$ denotes the first-order deformations
of $f$ modulo reparametrization --- the tangent space to
$\overline{M}_D(\com\proj^4/\com\proj^4_{\mathbb R},d)$ at $[f]$.
To
carry out the corresponding closed calculation, 
 the weights of each of the first two terms and computed and divided with
cancelling 0-weights.
In the open case, more information is needed about the
weight
0 components since
 isomorphisms of two copies of the trivial real
representation of $\mathbf{T}^n$ need not preserve orientation. 
We will use 
the real Euler sequence
to linearize the exact sequence \eqref{5t3} 
to get a better handle
on the sign.


Let $X = V_{\lambda_1/p}$ and 
$Y = \R \oplus V_{\lambda_1}\oplus V_{\lambda_3}.$
Let $$U = X\otimes \C, \ \
W = Y \otimes \C.$$
 Let $Y'
\subset Y$ denote the 2-dimensional linear subspace
corresponding to $V_{\lambda_1}$ and let $Y''$ denote the
$\mathbf{T}^n$-invariant complement. Consider the 
commutative diagram
in Figure \ref{fig:lin}.
\begin{figure}
\centering
\[
\xymatrix{ \rs(D,\O_{\C \PP^1}) \ar[r]\ar[d] & \rs(D,\O_{\C \PP^1}(1)
\otimes U) \ar[r]\ar[d]^{d\!\tilde f}
& \rs(D, T_{\C \PP^1}) \ar[d]^{df} \\
\rs(D,f^*\O_{\C \PP^4}) \ar[r]\ar[d] & \rs(D,f^*\O_{\C \PP^4}(1)\otimes W) \ar[r]\ar[d] & \rs(D,f^*T_{\C \PP^4}) \ar[d] \\
0 \ar[r] & \Def(D,f) \ar[r] & \Def(D,f) }
\]
\caption{} \label{fig:lin}
\end{figure}
All rows and columns are exact. The rows are the sections functor
applied to the Euler sequence. The rightmost column is obtained
from the deformation exact sequence by the identifications
\[
\Aut(D) = \rs(D, T_{\C \PP^1}), \qquad \Def(f) = \rs(D,f^*T_{\C
\PP^4}).
\]

In order to discuss the orientations of the vector spaces in
Figure \ref{fig:lin}, we must digress for a moment on the subject
of $Pin$ structures. In the following, $\p$ denotes a $Pin$
structure on $\C \proj^4_\R,$ $\p'$ denotes a $Pin$ structure on
$\O_{\C \PP^4}(5)_\R$ and $\s$ denotes a $Spin$ structure on $Q_\R.$
To define $N_d^{disk},$ a structure $\s$ was fixed.
By Lemma \ref{lm:sp} of Section \ref{sec:ecf}, 
any two of $\p,\p',\s,$ naturally determine
the third. As explained in the proof of Lemma \ref{lm:ro} of 
Section \ref{sec:ecf}, the
choice of $\p$ and $\p'$ compatible with $\s$ induces the
isomorphism $\det(\mathcal F_d) \simeq \L$ used to determine the
sign of $e(F_d).$ In particular, if we orient $N_f$ and $(F_d)_f$
using $\p$ and $\p'$ compatible with $\s,$ we will be calculating
the weights of the $\mathbf{T}^2$-action correctly.

In order to facilitate calculations, we choose $\p$ in a convenient
way. Then we let $\s$ and $\p$ induce $\p'.$ Indeed, choose first
a $Pin$ structure $\hat\p$ on $\O_{\C \proj^4}(1)_\R.$ Via the
natural homomorphism,
\[
Pin(1) \rarr Pin(5),
\]
induce a $Pin$ structure $\tilde \p$ on $\O_{\C \proj^4}(1)_\R
\otimes Y.$ In addition, equip $\O_{\C \PP^4_\R}$ with the canonical
$Spin$ structure. By \cite[Lemma 8.1]{jake1}, via
the Euler sequence, $\tilde \p$ induces a $Pin$ structure $\p$ on
$T_{\C \proj^4_\R}.$ Since we have chosen $Pin$ structures
compatibly, by \cite[Lemma 8.4]{jake1}, the middle row of the
diagram in Figure \ref{fig:lin} respects orientation. The columns
respect orientation by definition.

At this point, we focus attention on the middle column of Figure
\ref{fig:lin}, which is the desired linearization of the
deformation exact sequence. Because of the way we have induced
$\tilde \p$ from $\hat \p,$ by an argument similar to the proof of
\cite[Lemma 8.4]{jake1}, we may assume that the natural
isomorphism
\[
\rs(D,f^*\O_{\C \PP^4}(1)\otimes W) \simeq \rs(D,f^*\O_{\C
\PP^4}(1))\otimes Y
\]
preserves orientation. Here, we have used the
right-to-left lexicographical orientation of the tensor product.
Up to a degree independent sign, we may assume that the
isomorphism
\[
\rs(D,\O(1) \otimes U) \simeq \rs(D,\O(1)) \otimes X
\]
also preserves sign. So, it remains to compute the third term of
the exact sequence,
\begin{equation}\label{eq:lses}
0 \rarr \rs(D,\O(1))\otimes X \overset{d \!\tilde f}{\rarr}
\rs(D,\O(p))\otimes Y \rarr \Def(D,f) \rarr 0,
\end{equation}
from the first two.

We assume without loss of generality that the induced $\mathbf{T}^2$-action
on
\[
\rs(D,f^*\O_{\C \PP^4}(1)) \simeq \rs(D,\O(p))
\]
has weights
\[
\frac{\lambda_1}{p},\frac{3\lambda_1}{p},\ldots,
\]
as opposed to their negatives. That is, we assume the action on
the underlying vector space considered in Section \ref{sec:ws} has
weight $\lambda = \lambda_1/p.$ This depends on the choice of
$\hat \p.$ One way or another, the opposite disk has the opposite
sign, so we can always interchange $\zeta_1$ and $\zeta_2$ to
satisfy our assumption. Note, however, that the full localization
contribution of $[f]$ including the weights of the obstruction
bundle is invariant under the symmetry 
$$\lambda_1 \mapsto
\lambda_2 = -\lambda_1, \ \ \lambda_3 \mapsto \lambda_4=-\lambda_3.$$

Let $Y' \subset Y$ denote the 2-dimensional linear subspace
corresponding to $V_{\lambda_1}$ and let $Y'' \simeq V_{\lambda_3}
\oplus \R$ denote its $\mathbf{T}^2$-invariant complement. Note 
\[
\im(d \!\tilde f) \subset \rs(D,\O(p))\otimes Y'.
\]
We study the induced morphism,
\begin{equation*}
d \!\tilde f' : \rs(D,\O(1))\otimes X \rarr\rs(D,\O(p))\otimes Y',
\end{equation*}
carefully in order to calculate the weights of the cokernel with
attention to sign --- necessary because of the trivial
representations that occur in the domain and the range. Indeed, by
formulas \eqref{eq:t} and \eqref{eq:ws}, we have
\begin{gather}
\rs(D,\O(1))\otimes X \simeq V_0 \oplus V_{2\lambda_1/p} \label{eq:d}\\
\nonumber
\rs(D,\O(p))\otimes Y' \simeq \bigoplus_{i = 0}^{(p-1)/2} V_{(2i +
1 + p)\lambda_1/p } \oplus \bigoplus_{i = 0}^{(p-3)/2} V_{(2i + 1
- p)\lambda_1/p } \oplus V_0. 
\end{gather}
The map $d \!\tilde f'$ is determined up to homothety on the
summand $V_{2\lambda_1/p}$ by $\mathbf{T}^2$-equivariance. 
Since, homotheties of an even dimensional vector
space preserve orientation, we need calculate no further. However,
we need more information to determine $d \!\tilde f'$ on the
trivial representation summand $V_0.$  It is possible to explicitly
write down oriented bases of $X$ and $Y'$ and see that $d\!\tilde f'$
preserves orientation on the summand $V_0.$ However, that would lead to
considerable notational complications. Our strategy is to modify
the action of $\mathbf{T}^2$ on $X$ and $Y'$ so that $d \!\tilde f$ is
still equivariant, but there are {\em no} copies of the trivial
representation in the decomposition to irreducibles. To check
equivariance, we may work over the complex numbers, thus
simplifying formulas.

Choose coordinates $z_1$ and $z_2$ on $U$ such that under the
action of $\mathbf{T}^2$,
\[
z_1 \mapsto e^{i\lambda_1/p} z_1, \qquad z_2 \mapsto e^{-i\lambda_1/p}
z_2.
\]
Then $f$ is given explicitly by
\[
[z_1:z_2] \mapsto [z_1^p: z_2^p : 0 : 0 : 0 ].
\]
Let $e_1,e_2,$ be a basis of $U,$ dual to $z_1,z_2.$ A section
\[
\xi \in H^0(\C \PP^1,\O(1)\otimes U)
\]
takes the form,
\[
\xi = \xi_1 e_1 \oplus \xi_2 e_2,
\]
where $\xi_1,\xi_2,$ are linear functions on $U.$ Let $c_1,c_2,$
be the basis of $Y'$ corresponding to $e_1,e_2.$ Then
\[
d \!\tilde f' (\xi) = \xi_1 p z_1^{p-1} c_1 + \xi_2 p z_2^{p-1}
c_2.
\]
If we allow $\mathbf{T}^2$ to act on $X$ by $\lambda_1(1/p
+ \epsilon)$ and on $Y'$ by $\lambda_1(1 + \epsilon),$ then $d
\!\tilde f$ will remain equivariant. The summands $V_0$ in
decompositions \eqref{eq:d}  both change to
$V_{-\epsilon}.$ Choosing $\epsilon$ small enough, we may assume
no new trivial summands appear. Since the direct sum decomposition
doesn't change on the level of vector spaces as we change weights,
we conclude that $d \!\tilde f'$ maps $V_0$ to $V_0$ preserving
orientation. Hence, the cokernel of $d \!\tilde f'$ has
equivariant Euler class,
\begin{align*}
e(\coker(d \!\tilde f')) &= \lambda_1^{p-1} p^{1-p}
(p+1)(p+3)\ldots (2p)\times \\
& \qquad \qquad \qquad \times (1-p)(3-p)\ldots (-4) \\
& = (-1)^{(p-1)/2} 2^{p-1} p! p^{1-p} \lambda_1^{p-1}.
\end{align*}
Now, $\rs(D,\O(p))\otimes Y''$ contributes directly to $\Def(D,f)$
as follows. Decomposing $Y'' \simeq V_{\lambda_3} \oplus \R,$ and
using formulas \eqref{eq:t} and \eqref{eq:ws}, we calculate
\begin{align*}
e\left(\rs(D,\O(p))\otimes V_{\lambda_3}\right) &= \prod_{i =
0}^{(p-1)/2}\left(\left(1-\frac{2i}{p}\right)\lambda_1 -
\lambda_3\right) \times \\
& \qquad \qquad \qquad \times
\left(\left(1-\frac{2i}{p}\right)\lambda_1 + \lambda_3\right),
\end{align*}
and,
\[
e\left(\rs(D,\O(p))\otimes \R \right) = p!!
p^{(p+1)/2}\lambda_1^{(p+1)/2}.
\]
In conclusion,
\begin{align*}
e(N_f) & = (-1)^{(p-1)/2} 2^{p-1} p^{-(3p-1)/2} p! p!!
\lambda_1^{(3p-1)/2} \times \\
& \qquad \times \prod_{i =
0}^{(p-1)/2}\left(\left(1-\frac{2i}{p}\right)\lambda_1 -
\lambda_3\right)\left(\left(1-\frac{2i}{p}\right)\lambda_1 +
\lambda_3\right).
\end{align*}
On the other hand, by \eqref{eq:ws} the Euler class of the
obstruction bundle is just
\[
e\left((F_d)_f\right) = e(\rs(D,\O(5p)) = (5p)!! p^{-(5p+1)/2}
\lambda_1^{(5p+1)/2}.
\]
Here, we are assuming that the $Pin$ structure that $\p'$ induces
on
\[
\O(5p)_\R \simeq f^*\O_{\C \PP^4}(5)_\R
\]
agrees with $\p_{-1}.$ If for a given $\p'$ this is not true,
reversing $\p'$ changes the sign of $e\left((F_d)_f\right)$ by
$-1$ for all $p$ by \cite[Lemma 2.10]{jake1}.  
So, we reverse $\s,$ thus reversing $\p'.$ Hence, the dependence
of the total sign of $I(\zeta_1,p)$ on $\s$ as claimed. Combining
everything, we obtain
\[
I(\zeta_1,p) = 4 \frac{(-1)^{\frac{p-1}{2}}}{p}
\frac{\frac{(5p)!!}{p!p!!}\left(\frac{\lambda_1}{2p}\right)^{p+1}}{\prod_{i
= 0}^{(p-1)/2}\left(\left(1-\frac{2i}{p}\right)\lambda_1 -
\lambda_3\right)\left(\left(1-\frac{2i}{p}\right)\lambda_1 +
\lambda_3\right)}.
\]
The extra factor of $\frac{1}{p}$ comes from the orbifold
structure of the moduli space
$\widetilde{M}_D(\com\proj^4/\com\proj^4_{\mathbb
R},d)$ at $[f]$. \qed

\section{Proof of the Euler class formula}\label{sec:ecf}

\subsection{Construction
of \texorpdfstring{$\widetilde{M}_D(\com\proj^4/\com\proj^4_{\mathbb R},d)$}{M
tilde}}

A detailed construction was not required in \cite{jake1} to prove the
invariance of $N_d^{disk}$. However, to apply the Atiyah-Bott localization
formula as in the proof of Theorem  \ref{tm1} and Proposition \ref{der}, as well
as for the obstruction bundle argument in the proof of Theorem \ref{tm3}, we
need the following result.
\begin{pr}\label{pr:co}
$\widetilde{M}_D(\com\proj^4/\com\proj^4_{\mathbb R},d)$ is a
smooth closed orbifold and $F_d$ is a smooth orbibundle.
\end{pr}
\begin{proof}
First, we give a detailed definition of
$\widetilde{M}_D(\com\proj^4/\com\proj^4_{\mathbb R},d).$ The
moduli space $\overline{M}_D(\com\proj^4/\com\proj^4_{\mathbb
R},d)$ is a smooth orbifold with corners. A point in a corner of
codimension $k$ corresponds to an open stable map with $k$ disk
components. Such a stable map may have arbitrarily many
sphere components. However, since we consider open stable maps of
genus $0,$ each sphere component belongs to a tree of sphere
components attached to a unique disk. We define the \emph{total
degree} of a disk component to be its own degree plus the degree
of all attached spheres components.

We classify corners of
$\overline{M}_D(\com\proj^4/\com\proj^4_{\mathbb R},d)$ by the
intersection types $I$ of the open stable maps. The
intersection type $I = (T_I,\ell_I)$ of an open stable map $f$ of
genus zero consists of a tree $T_I$ and a labelling $\ell_I$ of
the vertices of $T_I$ by non-negative integers. The vertices of
the tree correspond to disk maps and the edges correspond to nodes
connecting two disks. The labelling of a vertex is the total
degree of the corresponding disk component. We denote by $|I|$ the
number of vertices of $T_I,$ ---  the codimension of the
corresponding corner. Let $I$ be an intersection type, and let $e$
be an edge of $T_I$ connecting vertices $v_1$ and $v_2.$ Gluing
$I$ at $e$ to obtain $I'$ means contracting $e$ so that $v_1$ and
$v_2$ become a single vertex $v$ of $I'$ and defining
$\ell_{I'}(v) = \ell_I(v_1) + \ell_I(v_2).$ We define a partial
ordering on the set of intersection types by $I_1 < I_2$ if $I_1$
may be obtained from $I_2$ by a sequence of gluings. We use the
notation $\partial_I{M}_D(\com\proj^4/\com\proj^4_{\mathbb R},d),$
to denote the corner of
$\overline{M}_D(\com\proj^4/\com\proj^4_{\mathbb R},d)$ of
intersection type $I.$

We refer to codimension one corners as boundary components. Points
of the boundary correspond to stable maps with two disk
components. The intersection type $I,\ |I| = 2,$ of such a stable
map is essentially an unordered pair $\{d_1,d_2\},$ where the
numbers $d_i$ are the total degrees of each disk component. Here
$d_1 + d_2 = d.$ For each possible such $I,$ we define a smooth
involution $c_I$ of the corresponding boundary component as
follows. Choose a vertex $v_I$ of $T_I.$ 
Since the total degree is odd, $v_I$ can be chosen as
the unique vertex of odd degree.
For $$[f] \in
\partial_I {M}_D(\com\proj^4/\com\proj^4_{\mathbb R},d),$$
define $c_I(f)$ to be the open stable map obtained by replacing
the disk component corresponding to $v_I$ and all attached sphere
maps by their conjugates.

We extend $c_I$ to the closure
\[
\overline{\partial_I M_D(\com\proj^4/\com\proj^4_{\mathbb R},d)} =
\bigcup_{I' \geq I} \partial_{I'}
M_D(\com\proj^4/\com\proj^4_{\mathbb R},d)
\]
as follows. Let $[f] \in \partial_{I'}
M_D(\com\proj^4/\com\proj^4_{\mathbb R},d)$ for some $I' > I.$
Let $\{v_{I'}^i\}$ be the set of vertices of $I'$ that glue to
form $v_I.$ We define $c_I(f)$ to be the open stable map obtained
by replacing the disk-maps corresponding to the $v_{I'}^i$ and all
attached spheres maps by their conjugates.

Next, we define an
equivalence relation $\sim$ on
$\overline{M}_D(\com\proj^4/\com\proj^4_{\mathbb R},d)$ as
follows. Let $[p],[q] \in \partial_{I'}
M_D(\com\proj^4/\com\proj^4_{\mathbb R},d).$ We write $p \sim q$
if there exists $I \leq I'$ with $|I| = 2$ such that $c_I(p) = q.$
Finally, we can define
\[
\widetilde{M}_D(\com\proj^4/\com\proj^4_{\mathbb R},d) =
\overline{M}_D(\com\proj^4/\com\proj^4_{\mathbb R},d)/\sim.
\]

We now prove
$\widetilde{M}_D(\com\proj^4/\com\proj^4_{\mathbb R},d)$ is a
closed orbifold. We refer the reader to \cite[Section 2]{fo} for a
quick review of orbifolds. We essentially follow the notation
established there. Let
\[
\pi : \overline{M}_D(\com\proj^4/\com\proj^4_{\mathbb R},d)\rarr
\widetilde{M}_D(\com\proj^4/\com\proj^4_{\mathbb R},d)
\]
denote the quotient projection. Being a closed orbifold is a local
property, so we restrict our attention to small neighborhood of a
point $$[p] \in \widetilde{M}_D(\com\proj^4/\com\proj^4_{\mathbb
R},d).$$ 
If $\pi^{-1}([p])$ lies away from the corners of
$\overline{M}_D(\com\proj^4/\com\proj^4_{\mathbb R},d),$ there is
nothing to prove.

If $\pi^{-1}([p])$ meets the corners,
we construct an
orbifold chart in a neighborhood of $[p]$ in a canonical way. Let
$[\tilde p] \in \pi^{-1}([p]).$ We assume 
 $$[\tilde p] \in
\partial_{I_p} M_D(\com\proj^4/\com\proj^4_{\mathbb R},d), \ \
k = |I_p| - 1 \geq 1.$$  
Let $(V_{\tilde p},\Gamma_{\tilde
p},\psi_{\tilde p})$ be an orbifold chart on
$\overline{M}_D(\com\proj^4/\com\proj^4_{\mathbb R},d)$ at $\tilde
p.$ Here $V_{\tilde p}$ is a neighborhood of $0$ in
\[
\R^n_{+k} = \{(x_1,\ldots,x_n) \in \R^n\, | \,x_i \geq 0, i =
1,\ldots k\},
\]
$\Gamma_{\tilde p}$ is a finite group acting on $V_{\tilde p}$,
 and
$\psi_{\tilde p}$ is a $\Gamma_{\tilde p}$ invariant homomorphism
from $V_{\tilde p}$ to a neighborhood of $[\tilde p].$ 
From the definition of the orbifold structure on
$\overline{M}_D(\com\proj^4/\com\proj^4_{\mathbb R},d),$ the
groups $\Gamma_{\tilde p}$ for $[\tilde p] \in \pi^{-1}([p])$ are all
isomorphic. So, we may define a group $\widetilde \Gamma_p$ with
isomorphisms $\widetilde \Gamma_p \simeq
\Gamma_{\tilde p}.$ We define,
\[
\partial_I V_{\tilde p} = \psi_{\tilde p}^{-1}\left(\partial_I M_D(\com\proj^4/\com\proj^4_{\mathbb
R},d)\right).
\]
By definition of an orbifold with corners,
$\partial_I V_{\tilde p}$ is contained in a subset of $\R^n_{+k}$
where $|I| - 1$ of the coordinates $x_1,\ldots,x_k,$ are zero. By
definition of a smooth map of an orbifold, possibly shrinking the
charts $V_{\tilde p},$  we may assume that $c_I$ induces a smooth
$\widetilde \Gamma_p$ equivariant involution of the disjoint union
of the corners $\partial_{I'} V_{\tilde p},\, I' \geq I$ over
all $[\tilde p] \in \pi^{-1}([p]).$

Let $I_1,\ldots, I_k$ enumerate the set of $I \leq I_p$ such that
$|I| = 2,$ the intersection types of the boundary components
adjacent to $p.$ For a multi-index
\[
E = (\epsilon_1,\ldots,\epsilon_k) \in (\Z/2\Z)^k,
\]
we define an involution $c_E$ of
$\partial_{I_p}M_D(\com\proj^4/\com\proj^4_{\mathbb R},d)$ by
\[
c_{E} = \prod_{j = 1}^k c_{I_j}^{\epsilon_j}.
\]
So, the group $(\Z/2\Z)^k$ acts on
$\partial_{I_p}M_D(\com\proj^4/\com\proj^4_{\mathbb R},d),$ and by
definition of $\pi,$ acts transitively on $\pi^{-1}([p]).$ Define
\[
G_p \subset (\Z/2\Z)^k
\]
to be the stationary subgroup of $[\tilde p] \in \pi^{-1}([p]).$ The
definition does not depend on the choice of $[\tilde p]$ because we
are considering a transitive action of an abelian group. Define
\[
\hat V_p = \coprod_{[\tilde p] \in \pi^{-1}([p])} V_{\tilde p} \times
G_{p}, \qquad V_p = \hat V_p /\sim^*,
\]
where the equivalence relation $\sim^*$ is defined as follows. Let
$$E_j \in (\Z/2\Z)^k$$ 
denote the multi-index with $\epsilon_i = 0$
for $i \neq j$ and $\epsilon_j = 1.$ For $$(q,E), (q',E') \in \hat
V_p$$ define $(q,E) \sim^* (q',E')$ if $q = c_{I_j}(q')$ and $E =
E' + E_j$ for some $j.$ Define
\[
\Gamma_p = \widetilde \Gamma_{p} \times G_p.
\]
The group $\Gamma_p$ acts naturally on $V_p.$ Finally, let
$\psi_p$ be the $\Gamma_p$ invariant map from $V_p$ to
$\widetilde{M}_D(\com\proj^4/\com\proj^4_{\mathbb R},d)$ naturally
induced by the $\psi_{\tilde p}.$ Since $V_p$ is obtained by
gluing together $2^k$ neighborhoods of $0 \in \R^n_{+k}$ on
matching corners, $V_p$ is a neighborhood of $0\in\R^n.$ The 
triple $(V_p,\Gamma_p,\psi_p)$ is easily seen to specify a
natural orbifold structure in a neighborhood of $p.$

The involutions $c_I$ lift naturally to the bundles
\[
\hat F_d|_{\overline{\partial_I M_D(\C \proj^4/\C \proj^4_\R,d)}}.
\]
Then, the exact same proof extends to construct the structure of
an orbibundle on $F_d.$
\end{proof}

In the following, $\p$ denotes a $Pin$ structure on $\C
\proj^4_\R,$ $\p'$ denotes a $Pin$ structure on $\O_{\C
\PP^4}(5)_\R$ and $\s$ denotes a $Spin$ structure on $Q_\R.$
\begin{lm}\label{lm:sp}
Any two of $\s,\p,\p'$ determines the third.
\end{lm}
\begin{proof}
By the adjunction formula, the normal bundle $N_Q$ of $Q$ in $\C
\proj^4$ satisfies $N_Q \simeq \O(5)|_Q.$ So, the Lemma follows
from the exact sequence,
\[
0 \rarr TQ_\R \rarr T\C \proj^4_\R \rarr N_{Q_\R} \rarr 0,
\]
and \cite[Lemma 8.1]{jake1}.
\end{proof}

We denote  the determinant of the tangent bundle of
$\widetilde{M}_D(\com\proj^4/\com\proj^4_{\mathbb R},d)$ by $\L.$
\begin{lm}\label{lm:ro}
There exists a  topological isomorphism 
$\det{F_d} \simeq \L.$ Moreover, such an
isomorphism is determined canonically up to homotopy by the choice
of a $Spin$ structure on $Q_\R.$
\end{lm}
\begin{proof}
Choose $\p$ and $\p'$ compatible with $\s.$ Then $\p$ induces an
orientation on $\overline{M}_D(\com\proj^4/\com\proj^4_{\mathbb
R},d)$ and $\p'$ induces an orientation on $\hat F_d.$ So, there
exists a unique up to homotopy isomorphism,
\[
\det\left(\overline{M}_D(\com\proj^4/\com\proj^4_{\mathbb
R},d)\right) \simeq \det(\hat F_d),
\]
preserving orientation. We must check the identifications
involved in the construction of
$\widetilde{M}_D(\com\proj^4/\com\proj^4_{\mathbb R},d)$ are
compatible with the above isomorphism:
the sign of $c_I$ on
$\partial_I{M}_D(\com\proj^4/\com\proj^4_{\mathbb R},d)$ must be the
same as the sign of $c_I$ on $\hat
F_d|_{\partial_I{M}_D(\com\proj^4/\com\proj^4_{\mathbb R},d)}.$
This follows easily from \cite[Lemma 2.12]{jake1}.
\end{proof}

\subsection{Kuranishi structures}\footnote{Here, we use the
virtual moduli cycle construction of \cite{FOOO}, 
based on Kuranishi structures.
In Sections \ref{kur1}-\ref{kur2}, 
we assume the coordinate transforms and Kuranishi maps
are smooth, as explained in \cite[Appendix A1.4]{FOOO}.  
We anticipate that with
the completion of the generalized Fredholm theory currently being introduced in
\cite{HWZ1,HWZ2}, the ideas of Section \ref{kur1}-\ref{kur2} 
will translate into a
proof of Theorem \ref{tm3} based on that framework.}
\label{kur1}

We denote by $\overline{M}_D(Q/Q_\R,d)$ the moduli space of open stable maps to
$(Q,Q_\R).$ In general, the space $\overline{M}_D(Q/Q_\R,d)$ is a compact
metrizable space. In the following proof, we use the theory of
Kuranishi structures with corners developed in \cite{fo,FOOO} to
define intersection theory on $\overline{M}_D(Q/Q_\R,d).$ The
Kuranishi structures used here were shown to exist in
\cite{fo,FOOO}. See \cite[Appendix A and Section 7]{jake1} for a
very brief summary of this theory, from which we take our notational
conventions. In the following, unless explicitly noted, all Kuranishi structures
are Kuranishi structures with corners.

We will need the following definition, which is similar to the
notion of an involution of a Kuranishi structure \cite[Definition 7.1]{jake1},
but without property (E1).
Suppose $(X,\K)$ and $(X',\K')$ are spaces with Kuranishi structure
\begin{align*}
\K &= (V_p,E_p,\Gamma_p,s_p,\psi_p,V_{pq},h_{pq},\varphi_{pq},\hat\varphi_{pq}), \\
\K' &=
(V_p',E_p',\Gamma_p',s_p',\psi_p',V_{pq}',h_{pq}',\varphi_{pq}',\hat\varphi_{pq}').
\end{align*}
Let $f: X \rightarrow X'$ be a continuous map.
\begin{dfn}
An \emph{extension} $\tilde f$ of $f$ to a map of spaces with Kuranishi
structure consists of $\Gamma_p$-equivariant maps 
$$f_p : V_p \rightarrow
V_{f(p)}', \ \ \hat f_p : E_p \rightarrow E_{f(p)}'$$ covering
$f_p$ such that
\begin{enumerate}
\renewcommand{\theenumi}{(M\arabic{enumi})}
\renewcommand{\labelenumi}{\theenumi}
\item
$s_{f(p)}' \circ f_p = \hat f_p \circ s_p.$
\item
$\psi_{f(p)}' \circ f_p|_{s_p^{-1}(0)} = f \circ \psi_p.$
\item
$f_q$ maps $V_{pq} \subset V_q$ to $V_{f(p)f(q)}' \subset
V_{f(q)}'.$
\item
$f_p\circ\varphi_{pq} = \varphi_{f(p)f(q)}'\circ f_q$ and
$\hat f_p\circ\hat\varphi_{pq} =
\hat\varphi_{f(p)f(q)}'\circ\hat f_q'.$
\end{enumerate}
\end{dfn}
Now, suppose that $(X,\K)$ has a tangent bundle given by
$\Phi_{pq}$ and $(X',\K')$ has a tangent bundle given by
$\Phi_{pq}'.$ We say that $\tilde f$ is smooth if
\begin{equation*}
\Phi_{f(p)f(q)}'\circ \hat f_q =  \hat f_p \circ \Phi_{pq}.
\end{equation*}
We say that $\tilde f$ is an embedding if $f_p$ are all embeddings and $\hat
f_p$ are all injective bundle maps. 
The preceding definition of
maps of spaces with Kuranishi structure is very rigid and not likely to make a
very nice category. A better definition for a general morphism of spaces with
Kuranishi structure would be something like the diagram of Figure \ref{fig:hat}.

We will also need the notion of vector bundles over a space with Kuranishi
structure and their Euler classes. Vector bundles over a space with Kuranishi
structure were constructed in a very general form in \cite[Section 5]{fo} so as
to include the tangent bundle to a Kuranishi structure. The bundles we use here
correspond to the special case in which $F_1$ of \cite{fo} is
taken to have rank $0$ everywhere. We let $(X,\K)$ denote a general space with
Kuranishi structure as above.
\begin{dfn}
A \emph{vector bundle} $F$ over $(X,\K)$ consists of
\begin{enumerate}
\renewcommand{\theenumi}{(\arabic{enumi})}
\renewcommand{\labelenumi}{\theenumi}
\item
For each $p \in X,$ a $\Gamma_p$-equivariant vector
bundle $F^p \rightarrow V_p.$
\item
For each $p \in X$ and $q \in \im \psi_p,$ an $h_{pq}$-equivariant vector
bundle isomorphism $\Phi^F_{pq} : F_q|_{V_{pq}} \rightarrow
F_p|_{\im(\varphi_{pq})}$ covering $\varphi_{pq}.$
\end{enumerate}
\end{dfn}
The Euler
class of a vector bundle over a space with Kuranishi structure
should determine a cohomology class in a cohomology theory
for Kuranishi spaces.
However, such a theory has not been developed. Since the Euler
classes considered here always have critical dimension, essentially
all the information in the Euler class is contained in a single
number: the integral of the Euler class over the fundamental
class. So, we focus on defining the integral of the Euler
class. 

Let $F$
denote a vector bundle over $(X,\K)$ of rank equal to the expected
dimension of $(X,\K).$ Let $\L_\K$ denote the determinant of the
tangent bundle of $(X,\K),$ i.e., the line bundle over $(X,\K)$
determined locally by the line bundles
\[
\det(TV_p) \otimes \det(E_p)^* \rarr V_p.
\]
Assume that an isomorphism
\begin{equation}\label{eq:LKF}
\det(\L_\K) \simeq \det(F)
\end{equation}
has been specified. Choose a transverse perturbation of the space
with Kuranishi structure $(X,\K).$ See \cite[Theorem A.4]{jake1}
for a brief review of notation. Choose multi-valued sections
$\xi_p'$ of $F_p'$ such that multi-valued section $\xi_p' +
s_{p,n}'$ of $E_p' \oplus F_p'$ is transverse. Let $\sigma$ denote
the 0-dimensional rational simplicial complex determined by the
vanishing set of $\xi_p' + s_{p,n}'.$ The orientation of
$\sigma$ is determined by the isomorphism \eqref{eq:LKF}. Let
$|\sigma|$ denote the rational weighted cardinality of $\sigma.$
\begin{dfn}
We define the \emph{Euler class} of $F$ by
\[
\int_{[X,\K]} e(F) = |\sigma|.
\]
\end{dfn}
Straightforward cobordism arguments show the definition does
not depend on the choice of section $\xi$ or the perturbation of
Kuranishi structure. See \cite[Sections 4 and 17]{fo}.

\subsection{Proof of Theorem \ref{tm3}}\label{kur2}
We continue to employ the notation of the proof
of Proposition \ref{pr:co}. Since $\overline{M}_D(Q/Q_\R,d)$
consists of open stable maps, we may define
\[
\partial_I {M}_D(Q/Q_\R,d) \subset \overline{M}_D(Q/Q_\R,d)
\]
to be the subspace consisting of all open stable maps of
intersection type $I.$ Similarly, we may define the involution
$c_I$ of the corner of intersection type $I'$ for $I' \geq I,$ and
the quotient space $\widetilde{M}_D(Q/Q_\R,d).$ Let $\K_Q$ be a
Kuranishi structure with corners on $\overline{M}_D(Q/Q_\R,d).$ By
the arguments of \cite[Section 7]{jake1}, the involutions $c_I$
extend smoothly to the Kuranishi structure
$\K_Q|_{\overline{\partial_I M_D(Q/Q_\R,d)}}.$ Recapitulating the
proof of Proposition \ref{pr:co}, $\K_Q$ induces a Kuranishi
structure without boundary $\widetilde \K_Q$ on
$\widetilde{M}_D(Q/Q_\R,d).$

A transverse perturbation of the space with Kuranishi structure
$(\widetilde{M}_D(Q/Q_\R,d), \widetilde \K_Q)$ defines a
simplicial complex consisting of a finite number of $0$-simplices
with rational weights. The weighted count of the  $0$-simplices   is
$N_d^{disk}.$ While the definition is not exactly the same as
the definition given in \cite{jake1}, the equivalence is
 not hard to verify.

\begin{figure}
\centering
\includegraphics[scale = .90]{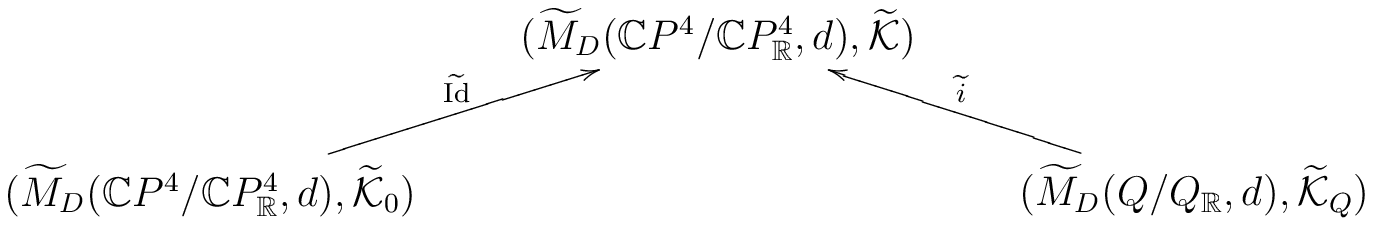}
\caption{}\label{fig:hat}
\end{figure}

An orbifold structure is a special case of a Kuranishi structure for
which all the bundles $E_p$ are rank 0. 
Let $\K_0$ denote the Kuranishi
structure on $\overline M_D(\C \proj^4/\C \proj^4_\R,d)$ 
coming from the orbifold
structure, 
and let $\widetilde \K_0$ denote the Kuranishi structure on
$\widetilde M_D(\C \proj^4/\C \proj^4_\R,d)$ coming from the 
orbifold structure
constructed in Proposition \ref{pr:co}. Also, let
\[
i : \widetilde{M}_D(Q/Q_\R,d) \rightarrow \widetilde M_D(\C \proj^4
/\C \proj^4_\R,d)
\]
denote the natural inclusion. We would like to construct a Kuranishi structure
$\widetilde \K$ on $\widetilde M_D(\C \PP^4/\C \PP^4_\R,d)$ for which 
 we have the
diagram of spaces with Kuranishi structure shown in Figure \ref{fig:hat}.

The structure $\widetilde \K$ may be obtained
from a Kuranishi structure on $\overline M_D(\C
\proj^4/\C \proj^4_\R,d)$ that admits extensions of the involutions $c_I$
for which the diagram of Figure \ref{fig:hat} holds with tildes
replaced by bars. For each $[p] \in \overline{M}_D(Q/Q_\R,d),$ we
extend the Kuranishi neighborhood
$$(V^Q_p,E^Q_p,\Gamma^Q_p,s^Q_p,\psi^Q_p)$$ given by $\K_Q$ to a
Kuranishi neighborhood $(V_p,E_p,\Gamma_p,s_p,\psi_p)$ for the space
$\overline M_D(\C \proj^4/\C \proj^4_\R,d).$ We detail the construction of
the extension for $p$ an irreducible stable map. The construction
for $p$ a reducible stable map is similar but notationally more
complicated.

Let $B_D(Q/Q_\R,d)$ and $B_D(\C \PP^4/\C \PP^4_\R,d)$
denote the Banach manifolds of $W^{1,p}$ maps
\[
(D,\partial D) \rightarrow (Q,Q_\R), \qquad (D,\partial D) \rightarrow (\C
\PP^4,\C \PP^4_\R),
\]
respectively. Let
\[
\E^Q \rightarrow B_D(Q/Q_\R,d), \qquad \E \rightarrow B_D(\C
\PP^4/\C \PP^4_\R,d),
\]
be defined fiberwise by
\begin{align*}
\E^Q_f &= L^p(D,f^* \Omega^{0,1}(TQ)), \quad f \in B_D(Q/Q_\R,d), \\
\E_f &= L^p(D,f^* \Omega^{0,1}(T\C P^4)), \quad f \in B_D(\C
\PP^4/\C \PP^4_\R,d).
\end{align*}
Let $\bar \partial^Q$ (resp. $\bar \partial$) denote the section
of $\E^Q$ (resp. $\E$) given by the non-linear Cauchy-Riemann
operator on maps to $Q$ (resp. $\C \PP^4$). Let $D\bar\partial^Q$
and $D\bar\partial$ denote choices of the vertical parts of their
respective linearizations.

We briefly outline the construction \cite{fo,FOOO} of the
Kuranishi neighborhood $(V^Q_p,E^Q_p,\Gamma^Q_p,s^Q_p,\psi^Q_p)$ in order to
explain how to extend it. First, choose a finite dimensional
subspace $\hat E^Q_p \subset \E^Q_p$ so that
\[
D\bar\partial^Q : T_p B_D(Q/Q_\R,d) \rightarrow \E^Q_p/\hat E^Q_p
\]
is surjective. Extend $\hat E^Q_p$ to a vector bundle $\check
E^Q_p$ over a neighborhood of $p$ in $B_D(Q/Q_\R,d)$ by parallel
transport. Over a sufficiently small neighborhood of $p,$
$D\bar\partial^Q$ surjects onto $\E^Q/\check E_p^Q.$ So, we may
define a smooth manifold
\[
\widetilde V^Q_p = (\bar\partial^Q)^{-1}(\check E_p^Q).
\]
We define $V^Q_p$ to be an appropriate section of the action of
the infinitesimal reparametrization group on $\widetilde V^Q_p.$
Then, we take $E^Q_p = \check E^Q_p|_{V^Q_p}$ and
\[
s_p^Q(x) = \bar\partial^Q(x) \in E^Q_p, \qquad x \in V^Q_p.
\]
The group $\Gamma_p^Q$ arises from the remaining discrete part of
the reparametrization group of $p.$ 

To extend the Kuranishi neighborhood
$(V^Q_p,E^Q_p,\Gamma^Q_p,s^Q_p,\psi^Q_p)$ to a Kuranishi neighborhood for
$\overline M_D(\C \PP^4/\C \PP^4_\R,d),$ we extend $\check E_p^Q$ by parallel
translation to a vector bundle $\check E_p$ over a neighborhood of
$p$ in $B_D(\C \PP^4/\C \PP^4_\R,d).$ Since $\C \PP^4$ is a homogeneous
space, the operator $D\bar\partial$ is surjective onto $\E_p$ at
every $[p] \in M_D(\C \PP^4/\C \PP^4_\R,d).$ So, we may define a smooth
manifold
\[
\widetilde V_p = (\bar\partial)^{-1}(\check E_p).
\]
The definitions of $V_p,\,E_p,\,s_p$ and $\Gamma_p,$ are just as
before. In the case when $p$ is a reducible stable map, the
extension procedure is similar except for one extra detail: we
must be careful to perform parallel translations using a
$c$-invariant metric so that the involutions $c_I$ extend to the
corners of the extended Kuranishi neighborhood. See \cite[Section
7]{jake1} for a discussion of this issue.

For each point in $[p] \in \overline{M}_D(Q/Q_\R,d)$ we have just
constructed a Kuranishi neighborhood
$(V_p,E_p,\Gamma_p,s_p,\psi_p)$ for $\overline M_D(\C \PP^4/\C
\PP^4_\R,d).$ It is straightforward to extend transition data
$V_{pq},\,h_{pq},\,\phi_{pq},\,\hat\phi_{pq},$ to these extended
neighborhoods. In order to complete the construction of $\K$ it
remains to define Kuranishi neighborhoods of points
\[
[p] \in \overline M_D(\C \PP^4/\C \PP^4_\R,d)\setminus \overline
M_D(Q/Q_\R,d).
\]
For such points, we take an orbifold coordinate chart as a
Kuranishi neighborhood, letting $E_p$ be a trivial bundle of rank
zero. Because the bundles $E_p$ are rank zero, it is again easy to
construct the associated transition data
$V_{pq},\,h_{pq},\,\phi_{pq},\,\hat\phi_{pq}.$ The existence of
embeddings as in the diagram of Figure \ref{fig:hat} follows
immediately from the construction of $\K.$

We proceed to extend the bundle $F_d$ of Section \ref{ssec:q3d} to
a bundle $\widetilde F_d$ over the space with Kuranishi structure
$(\widetilde M_D(\C \PP^4/\C \PP^4_\R,d),\widetilde\K).$ Again, it
suffices to construct a bundle $\check F_d$ over the space with
Kuranishi structure $(\overline M_D(\C P^4/\C P^4_\R,d),\K)$ so
that $\check F_d$ admits an extension of the involutions $c_I.$ We
define the fiber of $\check F_d^p$ at $f \in V_p$ as follows. By
construction, $f$ is a $W^{1,p}$-stable map satisfying the equation
$\bar\partial f \in \widetilde E_p.$ Since $f$ may not be
holomorphic, we have to be careful how we define the complex
structure on $f^* \O_{\C P^4}(5).$ Choose a $c$-invariant metric
on $\O_{\C P^4}(5).$ The associated complex connection on $\O_{\C
P^4}(5)$ will also be $c$-invariant. Equip $f^*\O_{\C P^4}(5)$
with the complex structure induced from the $(0,1)$ part of the
pull-back connection. Then, we define as before,
\[
\check F_d|_f = \rs(D,f^*\O_{\C P^4}(5)).
\]
The transition functions $\Phi^{F_d}_{pq}$ are tautological.
Because we used a $c$-invariant connection to induce the complex
structure on $f^*\O_{\C P^4}(5),$ the involutions $c_I$ lift to
involutions of $\check F_d|_{\partial_{I'}\overline M_D(\C P^4/\C
P^4_\R,d)}$ for $I' \geq I.$ Hence, $\check F_d$ descends to a
bundle $\widetilde F_d$ over $(\widetilde M_D(\C P^4/\C
P^4_\R,d),\widetilde \K).$ Furthermore, the proof of Lemma \ref{lm:ro} shows
that the isomorphism $\det(F_d) \simeq \L$ extends to an isomorphism
$\det(\widetilde F_d) \simeq \L_{\widetilde \K}.$

Let $\xi \in H^0(\C \PP^4, \O_{\C \PP^4}(5))$ be the section
defining the hypersurface $Q \subset \C \PP^4.$ For general $f \in V_p,$ the
pull-back $f^{-1}\xi$ is not a holomorphic section of $f^*\O_{\C
P^4}(5).$ Let
\[
P_f : W^{1,p}(D,f^*\O_{\C \PP^4}(5)) \rightarrow \rs(D,f^*\O_{\C
\PP^4}(5))
\]
denote the $L^2$ projection with respect to a $c$-invariant metric. We define a
section $\check \xi_p$ of $\check F_d^p$ by
\begin{equation}\label{eq:xif}
\check \xi_p (f) = P_f \left (f^{-1}\xi \right), \qquad f \in V_p.
\end{equation}
The local sections $\check \xi_p$ clearly match under the
transition functions $\Phi^{F_d}_{pq}$ to define a global section
$\check\xi$ of $\check F_d.$ Since $\xi$ is $c$-invariant and
$P_f$ is defined with respect to a $c$-invariant metric, we
conclude that $\check\xi$ is compatible with the involutions
$c_I.$ So, $\check \xi$ descends to a section $\tilde \xi$ of
$\widetilde F_d.$
\begin{lm}\label{lm:tr}
If the Kuranishi neighborhoods $V_p$ are chosen sufficiently small, then the
sections $\check\xi_p$ vanish precisely on the image of the embedding of spaces
with Kuranishi structure
\[
(\overline M_D(Q/Q_\R,d),\K^Q) \stackrel{\tilde i}{\rarr} (\overline M_D(\C
\PP^4/\C \PP^4_\R,d),\K).
\]
Moreover, the sections $\check\xi_p$ are transverse to zero.
\end{lm}
We postpone the proof of Lemma \ref{lm:tr} until we complete the
proof of Theorem \ref{tm3}. Indeed, we calculate
\[
\int_{[\widetilde M_D(\C \PP^4/\C \PP^4_\R,d),\widetilde
\K]}e(\widetilde F_d)
\]
in two different ways. First, we construct a transverse
perturbation $\mathfrak P_Q$ of $(\widetilde
M_D(Q/Q_\R,d),\widetilde \K^Q).$ Let $\sigma_Q$ denote the
$0$-dimensional simplicial complex defined by $\mathfrak P_Q.$ By
definition $N_d^{disk} = |\sigma_Q|.$ Now, extend $\mathfrak P_Q$
to a transverse perturbation $\mathfrak P$ of $(\widetilde M_D(\C
\PP^4/\C \PP^4_\R,d),\widetilde \K).$ By Lemma \ref{lm:tr}, the local
sections $\widetilde \xi_p' + \tilde s_{p,n}'$ of $\widetilde
F_d^p + \widetilde E_p$ vanish transversely exactly on $\sigma_Q.$
Assuming the following lemma, this implies that
\begin{equation}\label{eq:ec1}
\int_{[\widetilde M_D(\C \PP^4/\C \PP^4_\R,d),\widetilde
\K]}e(\widetilde F_d) = N_d^{disk},
\end{equation}
\begin{lm}\label{lm:ork}
The orientation induced on each zero simplex in $\sigma_Q$ by the Kuranishi
structure $\widetilde\K^Q$ agrees with the orientation induced by the Kuranishi
structure $\widetilde \K$ and the vector bundle $\widetilde F_d.$
\end{lm}
We postpone the proof of Lemma \ref{lm:ork} until the end of the section.
 
On the other hand, the sections $\widetilde s_p$ of $\widetilde
E_p$ are transverse without any perturbation since $\C \PP^4$ is
convex. So, we may choose a transverse multi-valued
section $\eta_0$ of $F_d,$ and extend it to a multi-valued section
$\tilde \eta$ of $\widetilde F_d$ over $(\widetilde M_D(\C \PP^4/\C
\PP^4_\R,d),\widetilde \K)$ such that the local sections $\tilde \eta_p +
\widetilde s_p$ will be transverse. This shows that
\begin{equation}\label{eq:ec2}
\int_{[\widetilde M_D(\C \PP^4/\C \PP^4_\R,d),\widetilde
\K]}e(\widetilde F_d) = \int_{{\widetilde M}_D(\C P^4/\C P^4_\R,d)}e(F_d).
\end{equation}
Combining equations \eqref{eq:ec1} and \eqref{eq:ec2}, we deduce
Theorem \ref{tm3}. \qed

\begin{proof}[Proof of Lemma \ref{lm:tr}]
We may focus on $\check\xi_p$ for $[p] \in \overline{M}_D(Q/Q_\R,d).$
Indeed, for other $p,$ the Kuranishi neighborhood $V_p$ is just an
orbifold chart on $\overline{M}_D(\C \PP^4/\C \PP^4_\R,d).$ So, for all $f \in
V_p$ the pull-back $f^{-1}\xi$ is holomorphic and non-zero. So,
$\check \xi_p(f)$ is never zero. In case 
$$[p] \in \overline{M}_D(Q/Q_\R,d),$$ we need to
show that $\check \xi_p$ vanishes transversely on
\[
V_p^Q \subset V_p,
\]
but nowhere else.

First, we establish some notation. Let $\nabla$ denote a complex $c$-invariant
connection on $\O_{\C \PP^4}(5)$ and let $\nabla^f$ denote its
pull-back to $f^*\O_{\C \PP^4}(5).$ Let $Y^p \rarr V_p$ and $Z^p \rarr V_p$
denote the Banach space bundles with fibers
\[
Y^p_f = W^{1,p}(D,f^*\O_{\C
\PP^4}(5)), \qquad Z^p_f = L^p\left(D,\Omega^{0,1}(f^*\O_{\C
\PP^4}(5))\right).
\]
Define a map of Banach-space bundles
\[
d'' : Y^p \rightarrow Z^p
\]
by
\[
d''_f = \left ( \nabla^f \right)^{0,1} : W^{1,p}(D,f^*\O_{\C
\PP^4}(5)) \rightarrow L^p\left(D,\Omega^{0,1}(f^*\O_{\C
\PP^4}(5))\right).
\]
Since $\ker(d''_f) = \rs(D,f^*\O_{\C \PP^4}(5)),$
and $d''_f$ is surjective for all $f,$ we have a short exact sequence
\[
0 \rarr \check F_d^p \rarr Y^p \stackrel{d''}{\rarr} Z^p \rarr 0
\]
Let 
\[
R : Z^p \rarr Y^p
\]
denote the unique right inverse of $d''$ such that the image of $R_f$ is the
$L^2$ complement of $\ker(d''_f).$ Let 
\[
P : Y^p \rarr \check F_d^p
\]
denote the $L^2$ projection. Define a section $\hat \xi_p$ of $Y^p$ by 
\[
\hat \xi_p(f) = f^{-1}\xi.
\]
Throughout the following, we use $\| \cdot \|$ to denote a context dependent
norm. That is, for sections of Banach space bundles, $\|\cdot \|$ is the
appropriate Banach space norm and for operators, $\| \cdot \|$ is the
appropriate operator norm.

Reformulating the definition of $\check\xi_p$ given in \eqref{eq:xif}, we have
\[
\check\xi_p = P \hat \xi_p = \hat \xi_p - R_f \circ d''_f (
\hat \xi_p).
\]
We will argue that for $f$ close enough to $p,$ we have 
\begin{equation}\label{eq:me}
\| R_f \circ d''_f (f^{-1}\xi ) \|_{1,p} \leq \epsilon \dist(f,V_p^Q)
\end{equation}
for arbitrary epsilon. Here, $\dist(\cdot,\cdot)$ denotes an arbitrary distance
function. Since $\xi$ vanishes transversely at $Q,$ we know that
$\hat \xi_p$ vanishes transversely on $V_p^Q.$
That is, we know that
\[
\| f^{-1}\xi \|_{1,p} \geq \epsilon_0 \dist(f,V_p^Q).
\]
Choosing $\epsilon < \epsilon_0,$ estimate \eqref{eq:me} shows that after we
perturb $\hat\xi_p$ by $R_f \circ d''_f (\hat\xi_p )$ to obtain $\check
\xi_p,$ it is still transverse and vanishes only on $V^Q_p.$ 

In order to prove estimate \eqref{eq:me}, we calculate $d''_f(f^{-1}\xi).$
Let $j$ denote the complex structure on $D,$ let $J$ denote the complex
structure on $\C \PP^4$ and let $I$ denote the complex structure on $\O_{\C
\PP^4}(5).$ Assume that
\begin{equation}\label{eq:bpf}
\bar \partial f = \eta_f \in  (\check E_p)_f.
\end{equation}
We calculate,
\begin{align}
d''_f(f^{-1}\xi) & = \nabla^f (f^{-1}\xi) + I\nabla^f
(f^{-1}\xi)\circ j \label{eq:calc}\\
& = \nabla \xi \circ du + I (\nabla \xi \circ du) \circ j \notag\\
& = \nabla \xi \circ du + I \left(\nabla \xi \circ (- J \circ du
\circ j + \eta_f)\right) \circ j \notag\\
& = I \nabla\xi \circ \eta_f \circ j. \notag
\end{align}
The third equality uses equation \eqref{eq:bpf} and the fourth
equality uses the holomorphicity of $\xi$ to cancel $I$ and $-J.$

Composition with $\nabla \xi$ defines a linear map
\[
\nabla \xi \circ : E_p \rightarrow Z^p.
\] 
Observe that for $f \in V_p^Q,$ we have
\[
(E_p)_f = (E_p^Q)_f \subset \E^Q_f = L^p\left(D,
\Omega^{0,1}(f^*TQ)\right).
\]
Since $\xi$ vanishes on $Q,$ we see that $\nabla \xi \circ$ maps
$(E_p)_f$ to zero for $f \in V_p^Q.$ By continuity, we have
\[
\| (\nabla \xi \circ)_f \| \leq C \dist(f,V_p^Q).
\]
So, we infer from calculation \eqref{eq:calc} that
\begin{equation}\label{eq:nxi}
\| d''_f(f^{-1}\xi) \| \leq C \dist(f,V_p^Q) \|\eta_f\|.
\end{equation}
Furthermore, we can assume a uniform bound 
\begin{equation}\label{eq:R}
\|R_f \| \leq C'
\end{equation}
for all $f \in V_p.$ For any $\epsilon'>0,$ we can choose $V_p$ so small that
for all $f\in V_p$ we have
\begin{equation}\label{eq:eta}
\| \eta_f \| \leq \epsilon'.
\end{equation} 
So, choosing $\epsilon'$ such that $\epsilon'C C' <  \epsilon$ and combining
estimates \eqref{eq:nxi}, \eqref{eq:R} and \eqref{eq:eta}, we conclude estimate
\eqref{eq:me}. This completes the proof of Lemma \ref{lm:tr}.
\end{proof}

\begin{proof}[Proof of Lemma \ref{lm:ork}]
The following proof is a generalization of the argument given in \cite[Section
8, Proposition 8.8]{jake1}. In the following, we will abbreviate
\[
TB := TB_D(\C P^4/\C P^4_\R,d), \qquad TB^Q := TB_D(Q/Q_\R,d).
\]
We continue to use the bundles $Y$ and $Z$ introduced in the proof of Lemma
\ref{lm:tr}.
Since $V_p^Q = \check \xi^{-1}(0) \subset V_p,$ and by Lemma \ref{lm:tr},
$\check \xi$ is transverse to $0,$ we conclude that $d\check \xi$ induces and
isomorphism
\[
\det(TV_p^Q) \stackrel{\sim}{\rightarrow} \det(TV_p)\otimes \det(\check F_d)^*.
\]
along $V_p^Q.$ Using the fact that by construction $E_p^Q = E_p|_{V_p^Q},$ we
can tensor the above isomorphism with $\det(E_p^Q)^*  = \det(E_p)^*$ to obtain
an isomorphism
\begin{equation}\label{eq:dxi}
\L_{\K^Q} \stackrel{d\xi}{\longrightarrow} \L_{\K}\otimes \det(\check F_d)^*.
\end{equation}
The lemma will follow if we show that isomorphism \eqref{eq:dxi} respects the
canonical orientations of each of the three deteminant bundles $\L_{\K^Q},
\L_{\K},$ and $\det(\check F_d).$ For this purpose, we introduce the
commutative diagram of vector bundles over $\tilde V^Q_p$ of Figure
\ref{fig:3d}.
\begin{figure}
\centering
\includegraphics[scale = .90]{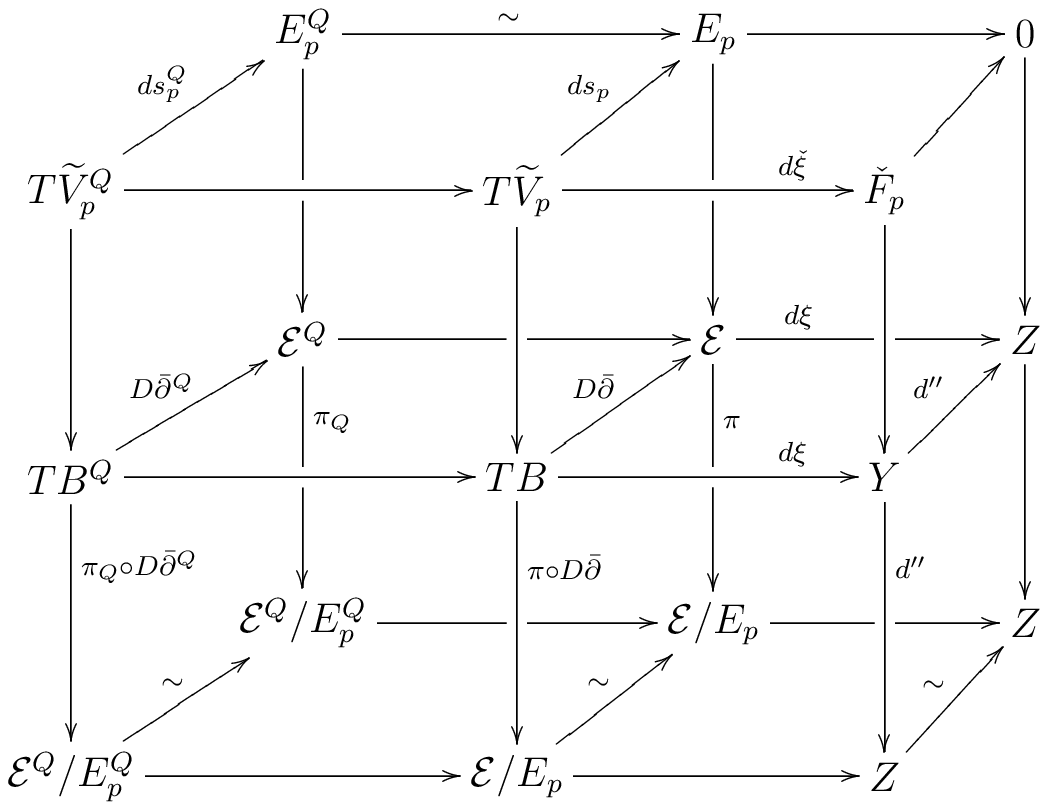}
\caption{}\label{fig:3d}
\end{figure}
Here, we implicitly consider the restrictions of $TB,\,TB^Q,\,\E,\,\E^Q,\,TV_p$
and $E_p$ to $V_p^Q.$ We denote by $\pi$ and $\pi_Q$ the canonical projections
to the quotient.

In order that all squares in the diagram of Figure \ref{fig:3d} commute, we must
choose the complex connection $\nabla$ on $\O_{\C P^4}(5)$ carefully. Indeed,
$d\xi|_Q$ induces an isomorphism of complex vector bundles
\[
d\xi|_Q : T\C P^4|_Q/TQ \stackrel{\sim}{\longrightarrow} \O_{\C P^4}(5).
\]
So, we are free to choose $\nabla$ to be the complex connection induced by
$d\xi|_Q$ from the connection on $T\C P^4|_Q/TQ$ induced by the Levi-Civita
connection of $\C P^4$ with respect to the standard Kahler metric. This ensures
that 
\[
d'' \circ d\xi =  d\xi \circ D\bar\partial,
\]
i.e. that the middle right horizontal square of the diagram commutes.
Moreover, the proof of Lemma \ref{lm:tr} implies that
\[
i_F \circ d\check\xi  =  d\xi \circ i,
\]
i.e. that the front upper right vertical square of the diagram commutes. The
commutativity of the remaining squares is straightforward. 

Moreover, all columns and rows of Figure \ref{fig:3d} are exact. We think of the
diagram as an exact square of two-step complexes. Applying the determinant
functor to the long exact sequence of a short exact sequence of two step
complexes yields an isomorphism of determinant bundles. In particular, the five
non-trivial rows and columns of Figure~\ref{fig:3d} give isomorphisms
\begin{gather}
\det(d s_p^Q) \simeq \det(D\bar\partial^Q), \quad \det(d s_p) \simeq
\det(D\bar\partial), \quad \det(\check F_d) \simeq \det(d''), \notag \\
\det(D\bar\partial) \simeq \det(D\bar\partial^Q)\otimes\det(d'')^*, \quad
\det(d s_p^Q) \simeq \det(d s_p) \otimes \det(\check F_d)^*.\label{eq:iso1}
\end{gather}
\begin{figure}
\centering
\includegraphics[scale = 1.00]{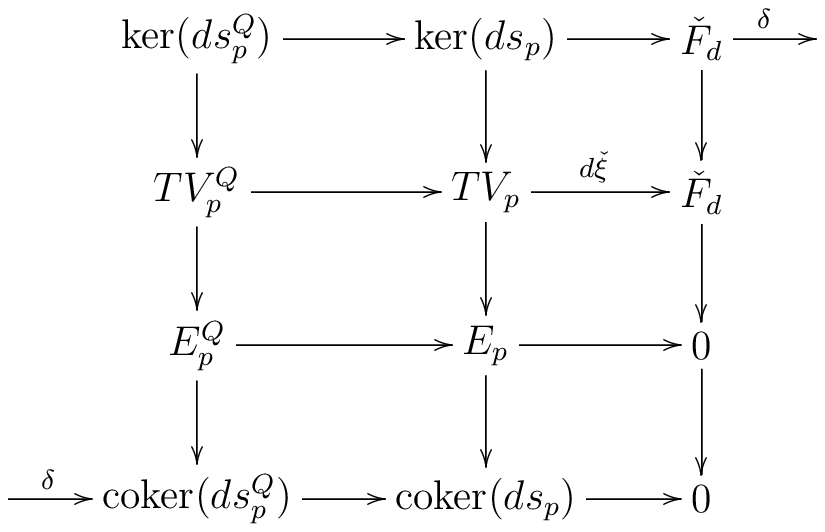}
\caption{}\label{fig:sn}
\end{figure}
Taking determinants of the rows and columns of the diagram of
Figure~\ref{fig:sn}, we obtain isomorphisms
\begin{gather}
\det(ds_p^Q) \simeq \L_{\K_Q}, \qquad \det(ds_p) \simeq \L_\K \notag\\
\L_{\K_Q} \simeq \L_\K \otimes \det(\check F_d)^*, \qquad \det(ds_p^Q) \simeq
\det(ds_p) \otimes \det(\check F_d)^*. \label{eq:iso2}
\end{gather}
Clearly the last isomorphism of \eqref{eq:iso1} agrees with the last
isomorphism of~\eqref{eq:iso2}. Putting all the isomorphisms of \eqref{eq:iso1}
and~\eqref{eq:iso2} together, we obtain the diagram of Figure \ref{fig:det}.
\begin{figure}
\centering
\[
\xymatrix{ \L_{\K_Q}  \ar[r]\ar[d]^{d\check\xi} &
\det(ds_p^Q)  \ar[r]\ar[d] & \det(D\bar\partial^Q)  \ar[d] \\
\L_\K \otimes \det(F_d)^* \ar[r] & \det(ds_p)\otimes \det(F_d)^* \ar[r] &
\det(D\bar\partial)\otimes \det(d'')}
\]
\caption{} \label{fig:det}
\end{figure}
The left square of Figure~\ref{fig:det} commutes by the commutativity of
Figure~\ref{fig:sn}. The right square of Figure~\ref{fig:det} commutes by the
commutativity of Figure~\ref{fig:3d}. The rows of Figure~\ref{fig:det} preserve
orientation by definition. The right column of Figure~\ref{fig:det}
preserves orientation by \cite[Section 8, Proposition 8.4]{jake1} and the
compatibility of $\s,\p,$ and $\p'.$ By commutativity the left column of
Figure~\ref{fig:det} also preserves orientation, completing the proof.
\end{proof}

\section{Multiple cover formula} \label{mcf}
\subsection{Local \texorpdfstring{$\com \proj^1$}{P1}}

Consider $\com\proj^1$ with the anti-holomorphic involution $c$
defined in Section \ref{ccc2}. The involution $c$ lifts canonically
 to
$\mathcal{O}_{\com\proj^1}(-1)$.
The total space of
the rank 2 bundle
$$\mathcal{O}_{\com\proj^1}(-1)
\oplus \mathcal{O}_{\com\proj^1}(-1) \rarr \com\proj^1$$
with the associated anti-holomorphic involution
may be viewed as local model for a rational curve in $Q$.

The local disk invariants of $\com\proj^1$ are, by definition,
$$L_d^{disk}= \int_{\widetilde{M}_D(\com\proj^1/\com\proj^1_{\mathbb{R}},d)}
e(G_d \oplus G_d),$$
where
$G_d$ is the real vector bundle over
with fiber
$$G_d|_{[f:(D,\partial D)\rarr(\com\proj^1,\com\proj^1_{\mathbb{R}})]}
= H^1(C,\tilde{f}^*\mathcal{O}_{\com\proj^1}(-1))_{\mathbb{R}}.$$
Here,
$$[\tilde{f}:C\rarr \com\proj^1] \in \overline{M}_{\mathbb R}(\com\proj^1,d)$$
is the stable rational map obtained from the stable
disk map via reflection.
As before, we consider only the $d$ odd case.

\begin{pr}\label{der}
For $d$ odd,
$L_d^{disk} ={2}{d^{-2}}.$
\end{pr}

The factor of 2 on the right occurs since the
original $\com\proj^1$ consists of 2 disks.
Hence, Proposition \ref{der} may be viewed as the calculation of
twice the multiple cover contribution of a single disk. 

We emphasize that $d$ is the degree of the stable rational map obtained from the
stable disk map by reflection. In particular, the degree of a map $f$ in
$\widetilde{M}_D(\com\proj^1/\com\proj^1_{\mathbb{R}},d)$ restricted to
$\partial D$ may be any odd integer less than or equal to $d.$ This is
necessary in order to perform the analogue of the construction of Proposition
\ref{pr:co} for target space $\C P^1,$ as the diffeomorphisms $c_I$ do not
preserve boundary degree. In this regard, our approach differs from the approach
of \cite{katzliu}.

\subsection{Torus action}

Let the torus $\mathbf{T}$
act on $\com\proj^1$
by
$$\xi \cdot [u,v] = [\xi u, \overline{\xi} v].$$
The fixed points are
$$\zeta_1=[1,0], \ \  \zeta_2=[0,1]$$
with
tangent weights
$2\lambda$ and $-2\lambda$
respectively.
The $\mathbf{T}$-action preserves
$\com\proj^1_{\mathbb{R}}$, and therefore determines
a translation action on the moduli space of disks
$\widetilde{M}_D(\com\proj^1/\com\proj^1_{\mathbb{R}},d)$.

The $\mathbf{T}$-action lifts to $\mathcal{O}_{\com\proj^1}(-1)$
with fiber
weights
$-\lambda$ and $\lambda$ over $\zeta_1$ and $\zeta_2$ respectively,
the unique lift which respects the
real structure on $\mathcal{O}_{\com\proj^1}(-1)$.

\subsection{Proof of Proposition \ref{der}} 
We give two different proofs.
\begin{proof}[First proof]
The invariants $L^{disk}_d$ are calculated by $\mathbf{T}$-equivariant
localization.
Since the steps are so similar to the proof of Theorem 1,
we give an abbreviated account.

As before, to each map $[f] \in
\overline{M}_D(\com\proj^1/\com\proj^1_{\mathbb{R}},d)^{\mathbf T}$,
we associate an intersection disk.
The intersection disk terms $I(\zeta_1,p)$ and
$I(\zeta_2,p)$ are both equal to
$$(-1)^{\frac{p-1}{2}} \frac{2^{1-e}}{p}\frac{(p!!)^2}{p!}.$$
The sign $(-1)^{\frac{p-1}{2}}$ comes from the normal bundle to the disk map
just as in the proof of Lemma \ref{ft1}. We do not have to calculate the
orientation of $H^1(C,\tilde{f}^*\mathcal{O}_{\com\proj^1}(-1))_{\mathbb{R}}$
because it appears twice.

Givental's equivariant correlator $S_L$ for the local
geometry is defined by
\begin{multline*}
S_L(T,\hb)=\\
\sum_{r\geq 0} e^{(H/\hbar+r) T}   e_{2*}(
\frac {c_{\text{top}}( H^1(\mathcal{O}(-1))\oplus
H^1(\mathcal{O}(-1)))} {\hb-\psi_2})\in H^*_{\mathbf{T}}(\com\proj^1)
\end{multline*}
where
$$e_2: \overline{M}_{0,2}(\com\proj^1,d) \rarr \com\proj^1$$
is the evaluation map.
The calculation of $S_L$ in \cite{g2} is much easier than $S_Q$,
$$S_L(T,\hb) = \sum_{r\geq 0} e^{rT}\frac{\prod_{s=0}^{r-1}
(H+s \hb)^2}{\prod_{s=1}^r (H-\lambda+s \hb)(H+\lambda+s\hb)}.$$
No mirror transform is needed.

The local disk potential $\mathcal{F}^{disk}_L$ is defined by
summing over odd degrees,
 $$\mathcal{F}^{disk}_L = \sum_{d\ odd} e^{dT/2} L^{disk}_d.$$
As before,
\begin{multline*}
\mathcal{F}^{disk}_L = \\
\sum_{p\ odd} \langle S_L(T,\frac{2}{p}\lambda), [\zeta_1] \rangle
\cdot I(\zeta_1,p)+
\sum_{p\ odd}\langle S_L(T,-\frac{2}{p}\lambda), [\zeta_2] \rangle
\cdot I(\zeta_2,p).
\end{multline*}
Evaluation yields
$$\mathcal{F}^{disk}_L = 2\sum_{r\geq 0} \sum_{p \ odd}
e^{(\frac{p}{2}+r)T}
\frac{2^{1-p-2r}}{(p+2r)^2} \frac{(-1)^{\frac{p-1}{2}}}{p}
 \frac{(p!!)^2}{p! r!} \frac{\prod_{i=1}^r (p+2i)^2}
{\prod_{i=1}^r (p+i)}.$$
The proof of the Proposition is concluded by extracting the
$e^{dT/2}$ terms on the right and executing the sum,
\begin{eqnarray*}
L_d^{disk} &=& \frac{2}{d^2}\sum_{1\leq p \ odd \leq d}
 \frac{(d!!)^2}{2^{d-1}}
\frac{(-1)^{\frac{p-1}{2}}}
{\left(\frac{d-p}{2}\right)!\left(\frac{d-p}{2}+p\right)!p }
\\
& = & \frac{2}{d^2}.
\end{eqnarray*}
Remarkably, the binomial identity required is
exactly \eqref{cccz}!
\end{proof}

\begin{proof}[Second proof]
Although each copy of $\O_{\C \PP^1}(-1)$ admits a unique real 
$\mathbf{T}$-action,
there is some freedom in the choice of action on the sum 
\[\O_{\C \PP^1}(-1)\oplus\O_{\C \PP^1}(-1).\]
Indeed, we can choose the action induced by the 
isomorphism of vector bundles
with real structure,
\begin{equation}\label{r5t}
\O_{\C \PP^1}(-1)\oplus\O_{\C \PP^1}(-1) \simeq \O_{\C \PP^1}(-1)
\otimes_\C W_\lambda.
\end{equation}
Here, $W_\lambda = V_\lambda \otimes \C$ and $V_\lambda$ is the real 
representation of the torus of weight $\lambda.$

With the $\mathbf{T}$-action determined
by the right side of \eqref{r5t}, the
localization contribution from any {\em reducible}
 open stable map vanishes. Indeed,
we work this out explicitly for a torus fixed open stable map 
$f$ with reducible
domain $D$ consisting of one disk component $D_o$ and one sphere component
$D_c.$ The general case is similar. 

Denote the single node of $D$ by $z.$ Consider the exact sequence
\begin{align*}
& 0 \rightarrow \rs(z,f^*\O_{\C \PP^1}(-1)\otimes W_\lambda) \rightarrow
\rh(D,f^*\O_{\C \PP^1}(-1)\otimes W_\lambda) \\
& \quad \rarr \rh(D_o,f^*\O_{\C \PP^1}(-1)\otimes W_\lambda)
\oplus H^1(D_c,f^*\O_{\C \PP^1}(-1)\otimes W_\lambda) \rarr 0.
\end{align*}
We claim that $\rs(z,f^*\O_{\C \PP^1}(-1)\otimes W_\lambda)$ contains a zero
weight space and therefore so does $\rh(D,f^*\O_{\C \PP^1}
(-1)\otimes W_\lambda).$
which
immediately implies that the localization contribution of $[f]$ vanishes.
Indeed, $$\rs(z,f^*\O_{\C P^1}(-1)\otimes W_\lambda)$$ 
is just the fiber of
$f^*\O_{\C P^1}(-1)\otimes W_\lambda$ at $z.$ Denoting by $\C_\lambda$ the
complex representation of the torus of weight $\lambda$ we have 
\[
W_\lambda \simeq \C_\lambda \oplus \C_{-\lambda}, \qquad  f^*\O_{\C P^1}(-1)_z
\simeq \C_{\pm \lambda}.
\]
Hence,
\[
\rs(z,f^*\O_{\C P^1}(-1)\otimes W_\lambda) = \C_0 \oplus \C_{\pm 2\lambda}.
\]
as claimed.

It remains to calculate the localization contribution from the single torus
fixed irreducible open stable map of degree $d.$ This is easily seen to be
$\frac{2}{d^2}.$
\end{proof}

The two proofs of Proposition \ref{der} together provide a geometric
evaluation of the binomial sum of Lemma \ref{idd}.

\subsection{Integrality}
The virtual disk counts $n^{disk}_d$
of Definition \ref{tm2} have not yet been proven to be integers. However,
the table in Section \ref{xxx} provides substantial evidence
for the integrality claim.


\section{Tables}
\label{xxx} Table \ref{t1} shows the value of the disk
Gromov-Witten invariant $N_d^{disk}$ for small $d.$

\begin{table}[h]
\caption{}
\label{t1}
\begin{center}
\begin{tabular}{c|l}
$d$& $N_d^{disk}$ \\
\hline
1& $30$ \\[.2cm]
3&$\frac{4600}{3}$ \\[.2cm]
5&$\frac{5441256}{5}$ \\[.2cm]
7&$\frac{47823842250}{49} $\\[.2cm]
9&$\frac{28973369597500}{27} $\\[.2cm]
11&$\frac{160812279574853640}{121}$ \\[.2cm]
13&$\frac{301152359429255569200}{169} $\\[.2cm]
15&$2528247216911976710478 $ \\[.2cm]
17&$\frac{1081454384062665012504422250}{289}$ \\[.2cm]
19&$\frac{2066166201384849550431238897500}{361}$ \\[.2cm]
21&$\frac{440336544802747748968402664543390}{49}$ \\[.2cm]
23&$\frac{7625558614788648016004683159051585650}{529}$ \\[.2cm]
25&$\frac{2942308498496733293257158606365620128756}{125}$ \\[.2cm]
27&$\frac{9481608375404186315963625791852891724001750}{243}$ \\[.2cm]
29&$\frac{55101515400393595065761084565358564820821590000}{841}$
\end{tabular}
\end{center}
\end{table}

\clearpage

Table \ref{t2}
shows the corresponding virtually enumerative invariants
$n_d^{disk}.$ Recall the virtual counts $n_d^{real}$ 
of real curves in
$Q$ differ by a factor of $1/2$,
$$n_d^{real}= \frac{1}{2} n_d^{disk}.$$

\begin{table}[h]
\caption{}
\label{t2}
\begin{center}
\begin{tabular}{c|l}
$d$ & $n_d^{disk}$ \\
\hline
1&30\\
3& 1530\\
5& 1088250\\
7& 975996780\\
9& 1073087762700\\
11& 1329027103924410\\
13& 1781966623841748930\\
15& 2528247216911976589500\\
17& 3742056692258356444651980\\
19& 5723452081398475208950800270\\
21& 8986460098015260183028517362890\\
23& 14415044640432226873354788580437780\\
25& 23538467987973866346057268850924917500\\
27& 39018964507836157678862657579522297754750\\
29& 65519043282275380577599387116954298241167170
\end{tabular}
\end{center}
\end{table}

\vspace{10pt}
\noindent Department of Mathematics \\
\noindent Princeton University \\
\noindent Princeton, NJ 08544\\
\noindent rahulp@math.princeton.edu

\vspace{10pt}
\noindent School of Mathematics \\
\noindent Institute for Advanced Study \\
\noindent Princeton, NJ 08540 \\
\noindent jake@ias.edu

\vspace{10pt}
\noindent School of Natural Science \\
\noindent Institute for Advanced Study \\
\noindent Princeton, NJ 08540\\
\noindent walcher@ias.edu


\begin{thebibliography}{[COGP]}

\bibitem{am} P. Aspinwall and D. Morrison, {\em Topological
field theory and rational curves}, Comm. Math. Phys. {\bf 151}
(1993), 245-262.

\bibitem{cogp} P. Candelas, X.  de la Ossa, P.  Green and L.  Parkes, {\em A
pair of Calabi-Yau manifolds as an exactly soluble superconformal field theory},
Nuclear Physics {\bf B359} (1991), 21-74.

\bibitem{FOOO} K. Fukaya, Y.-G. Oh, H. Ohto, K. Ono, {\em Lagrangian
intersection Floer theory, anomaly and obstruction}, Kyoto
University preprint, 2006.

\bibitem{fo}{K. Fukaya, K. Ono,}
    {\it Arnold conjecture and Gromov-Witten invariant,}
    {Topology {\bf 38} (1999), no. 5, 933--1048.}

\bibitem{g1} A. Givental, {\em Equivariant Gromov-Witten invariants},
Int. Math. Res. Notices {\bf 13} (1996), 613-663.


\bibitem{g2} A. Givental, {\em Elliptic Gromov-Witten invariants
and the generalized mirror conjecture}, math.AG/9803053.

\bibitem{grzas} T. Graber and E. Zaslow, {\em Open string
Gromov-Witten theory: calculation and a mirror theorem},
hep-th/0109075.

\bibitem{HWZ1} H. Hofer, K. Wysocki, E. Zehnder, {\em A General Fredholm Theory
I: A Splicing-Based Differential Geometry,} to appear in the Journal of the
European Mathematical Society, arXiv:math/0612604.

\bibitem{HWZ2} H. Hofer, K. Wysocki, E. Zehnder, {\em A General Fredholm Theory
II: Implicit Function Theorems,} arXiv:0705.1310.

\bibitem{katzliu} S. Katz and M. Liu, {\em Enumerative geometry
of stable maps with Lagrangian boundary conditions and multiple
covers of the disk}, Adv. Theor. Math. Phys. {\bf 5} (2002), 1-49.

\bibitem{ko}  M. Kontsevich, {\em Enumeration of rational curves via
torus actions}, in {\em The moduli space of curves}, (R.
Dijkgraaf, C. Faber, and G. van der Geer, eds.), Birkhauser, 1995,
335-368.

\bibitem{hams}{M. Kontsevich,}
    {\it Homological algebra of mirror symmetry.}
    {Proceedings of the International Congress of Mathematicians, Vol. 1, 2
    (Zurich, 1994), 120--139, Birkhauser, Basel, 1995.}

\bibitem{lly} B. Lian, K. Liu, and S.-T. Yau, {\em Mirror principle
I}, Asian J. Math. {\bf{4}} (1997), 729-763.

\bibitem{melissa} M. Liu, {\em Moduli of J-holomorphic curves
with Lagrangian boundary conditions and open Gromov-Witten invariants
for a $S^1$-equivariant pair}, math.SG/0210257.

\bibitem{oovafa} H. Ooguri and C. Vafa, {\em Knot invariants
and topological strings}, Nucl. Phys. B {\bf 577} (2000) 419-438. 

\bibitem{pan} R. Pandharipande, {\em Rational curves on hypersurfaces [after
A. Givental]}, S\'eminaire Bourbaki, 50\`eme ann\'ee, 1997-1998, no. 848.

\bibitem{RT}{Y. Ruan, G. Tian,}
    {\it A mathematical theory of quantum cohomology,}
    {J. Differential Geom. {\bf 42} (1995), no. 2, 259--367.}

\bibitem{S}{P. Seidel,} personal communication based on a remark of D. Joyce
and a talk of K. Fukaya at Northwestern in spring 2004.   
    
\bibitem{jake1} J. Solomon, {\em Intersection theory on the moduli
space of holomorphic curves with Lagrangian boundary conditions},
math.SG/0606429.

\bibitem{Walcher} J. Walcher, {\em Opening mirror symmetry on the
quintic}, hep-th/0605162.

\bibitem{We1}{J.-Y. Welschinger,}
    {\it Invariants of real symplectic 4-manifolds and lower bounds in real enumerative geometry,}
    {Invent. Math. {\bf 162} (2005), 195 - 234.}
   

\bibitem{We2}{ J.-Y. Welschinger,}
    {\it Spinor states of real rational curves in real algebraic convex 3-manifolds and enumerative invariants,}
        {Duke Math. J.  {\bf 127}  (2005), 89--121.}
        
\bibitem{WW}{K. Wehrheim, C. Woodward,}{\it Orientations for pseudo-holomorphic
	quilts,}
	{preprint.}

\bibitem{W}{E. Witten,}
    {\it Chern-Simons gauge theory as a string theory,}
    {The Floer memorial volume, 637--678, Progr. Math., {\bf 133},
    Birkhauser, Basel, 1995.}
    {arXiv: hep-th/9207094.}
\end{thebibliography}
\end{document}